\documentclass[12pt,reqno]{amsart}
\usepackage{mathrsfs, amsmath,amsfonts,amssymb}
\usepackage{amsmath}
\usepackage{stmaryrd}
\usepackage{mathtools}
\usepackage{fullpage}
\usepackage{color}
\allowdisplaybreaks
\newtheorem{proposition}{Proposition}
\newtheorem{corollary}{Corollary}

\newtheorem{theorem}{Theorem}
\newtheorem{remark}{Remark}
\usepackage{enumitem}
\usepackage{adjustbox}
\usepackage{hyperref}
\setlist{topsep=1.0em, itemsep=0.5em}

\newcommand{\ceil}[1]{\lceil {#1} \rceil}

\makeatletter
\renewcommand{\tocsection}[3]{%
	\indentlabel{\@ifnotempty{#2}{\bfseries\ignorespaces#1 #2.\,\,}}\bfseries#3}
\renewcommand{\tocsubsection}[3]{%
	\indentlabel{\@ifnotempty{#2}{\ignorespaces#1 #2\quad}}#3}
\renewcommand{\tocsubsubsection}[3]{%
	\quad\quad\quad\indentlabel{\@ifnotempty{#2}{\ignorespaces#1 #2\quad}}#3}

\newcommand\@dotsep{4.5}
\def\@tocline#1#2#3#4#5#6#7{\relax
	\ifnum #1>\c@tocdepth 
	\else
		\par \addpenalty\@secpenalty\addvspace{#2}%
		\begingroup \hyphenpenalty\@M
		\@ifempty{#4}{%
			\@tempdima\csname r@tocindent\number#1\endcsname\relax
		}{%
			\@tempdima#4\relax
		}%
		\parindent\z@ \leftskip#3\relax \advance\leftskip\@tempdima\relax
		\rightskip\@pnumwidth plus1em \parfillskip-\@pnumwidth
		#5\leavevmode\hskip-\@tempdima{#6}\nobreak
		\leaders\hbox{$\m@th\mkern \@dotsep mu\hbox{.}\mkern \@dotsep mu$}\hfill
		\nobreak
		\hbox to\@pnumwidth{\@tocpagenum{\ifnum#1=1\bfseries\fi#7}}\par
		\nobreak
		\endgroup
	\fi}
\AtBeginDocument{%
	\expandafter\renewcommand\csname r@tocindent0\endcsname{0pt}
}

\def\l@subsection{\@tocline{2}{0pt}{2.5pc}{5pc}{}}
\makeatother

\makeatletter
\newcommand{\oset}[3][0.25ex]{%
  \mathrel{\mathop{#3}\limits^{
    \vbox to#1{\kern-2\ex@
    \hbox{$\scriptstyle#2$}\vss}}}}
\makeatother
\title[]{finite free probability and $S$ transforms \\ of compact Jacobi processes}

\author[Nizar Demni]{Nizar Demni ${}^{*}$}
\address{${}^{*}$Division of Science and Mathematics, New York University Abu Dhabi, P.O. Box 129188, Abu Dhabi, United
Arab Emirate.s}
\email{nd2889@nyu.edu}

\author{Nicolas Gilliers ${}^{*,\,\dagger}$}
\address{${}^{\dagger}$ Université Paris Cité, CNRS, MAP5, F-75006 Paris, France}
\email{nicolas.gilliers@u-paris.fr} 

\author[]{Tarek Hamdi${}^{\triangle}$}
\address{${}^{\triangle}$Department of Management Information Systems  \\ College of Business and Economics\\ Qassim University \\ Saudi Arabia
and Laboratoire d'Analyse Math\'ematiques et applications LR11ES11 \\ Universit\'e de Tunis El-Manar \\ Tunisie\\}

\email{t.hamdi@qu.edu.sa} 


\usepackage[backend=biber,style=alphabetic,sorting=nyt]{biblatex}
\begin{document}
\maketitle
\begin{abstract}
We calculate the averaged characteristic polynomial and its finite $S$ transform for the Hermitian Jacobi process at any fixed time $t$. We give a direct proof that this sequence of polynomials solves the backward heat equation linked to the one-dimensional Jacobi operator. We also expand the averaged characteristic polynomials in terms of Jacobi polynomials, using the dual Cauchy identity for multivariate Jacobi polynomials and their mutual orthogonality. 

The finite free $S$ transform, introduced recently in \cite{arizmendi2024s}, is the finite free version of the free $S$ transform in that it behaves the same way with respect to the (finite) free multiplicative convolution. We present a finite difference and differential equation that the finite free $S$ transform of the averaged characteristic polynomials of the Hermitian Jacobi Process satisfies. In the high-dimensional limit, this yields a differential equation for the free $S$- transform of the free Jacobi process. We also prove a general technical lemma about the convergence of the finite differences of the finite free $ S$- transform.

\end{abstract}

\tableofcontents

\section{Introduction and main results}
 Inspired by free probability theory, Marcus, Spielman and Srivastava introduced finite free probability (hereafter FFP) theory and used it to solve the Kadison-Singer conjecture and to construct bipartite Ramanujan graphs (see \cite{MR4408504} and references therein). At the heart of this theory lie, among other things, convolution operations defined by Walsh and Szegő on real-rooted polynomials (see \cite{MR4408504} and references therein). When the latter are averaged characteristic polynomials of complex random matrices whose laws enjoy prescribed symmetry invariances, connections with Voiculescu’s free additive and multiplicative convolutions are made by considering the counting measure of the roots: finite free convolutions approximate, as the degree of the polynomials increases to infinity, the free convolutions. 

FFP theory profits from these connections, as much progress comes from finite analogues of results and transforms in the realm of free probability theory. 
For now, we highlight one example that will be connected to one of our forthcoming results. The finite free unitary convolution was introduced in \cite{MR4285332} where 
the author also derives the corresponding central limit theorem. The limiting polynomial is the so-called \emph{unitary Hermite polynomial}, which was shown in \cite{MR4912666} to coincide (up to rescaling) with the characteristic polynomial of Brownian motion on the unitary groups. The empirical root distribution converges weakly to the spectral distribution of the free unitary Brownian motion \cite{MR4912666}. 
We found it relevant to mention at this stage this polynomial since we shall prove in Proposition \ref{prop:relationunitaryjacobi} below a connection, via the Szegö transform, between it and the characteristic polynomial of the Hermitian Jacobi process. This connection is the finite free counterpart of a result of the same nature in the realm of free probability as stated in \cite{MR3071702}, Corollary 3.3.

$\beta$-diffusions are interacting particle systems that extrapolate eigenvalues of matrix-valued processes through a parameter $\beta > 0$ which may be interpreted as the inverse temperature in the Coulomb gaz picture (\cite{assiotis2023exact}, \cite{For}). It is sometimes called the Dyson index while $2/\beta$ is referred to as the Jack parameter in algebraic combinatorics (\cite{For}). Examples of $\beta$-diffusions include the so-called $\beta$-Pearson class which extends the famous Pearson class to higher dimensions. 
In particular, the $\beta$-Jacobi process was introduced and studied in \cite{MR2719370} in analogy with the already existing $\beta$-Dyson and $\beta$-Laguerre processes. 
Indeed, this particle system is an extrapolation to all $\beta > 0$ of the eigenvalues of the real symmetric ($\beta = 1$) and of the Hermitian Jacobi processes ($\beta=2$). 
A common property shared by these $\beta$-diffusions, as recently observed in \cite{auer2026hua} for the $\beta$-Hua-Pickrell diffusion, is the independence of their averaged characteristic polynomials of $\beta$ (\cite{MR4673381}, \cite{MR4781073}. Actually, the averaged elementary symmetric functions in these $\beta$-diffusions coincide with the elementary symmetric functions of the deterministic particle systems arising from the SDEs they satisfy
(informally and after rescaling) in the freezing regime $\beta \rightarrow +\infty$. As far as the $\beta$-Dyson diffusion is concerned, the corresponding deterministic particle system is a particular vortex system on the real line which coincides with the zeroes of a polynomial evolving with respect to the heat 
flow. This last fact was noticed for the first time by T. Tao\footnote{https://terrytao.wordpress.com/2017/10/17/heat-flow-and-zeroes-of-polynomials} and extended to the $\beta$-Pearson class 
in \cite{MR4673381}, \cite{MR4781073}, \cite{auer2026hua}. Quite amazingly, it was proved in \cite{FFPP} that the deterministic dynamics arising from freezing the $\beta$-Dyson, the $\beta$-Laguerre and the $\beta$-Jacobi processes also govern the evolution of the zeroes of averaged characteristic polynomials of certain random matrix models. 
Note however that the Jacobi process considered in \cite{FFPP}, section 5.1 is different from the Hermitian Jacobi process studied in \cite{MR3842169}, \cite{MR4173022}. 


In this paper, we study the high-dimensional behaviour of the roots of the averaged characteristic polynomial of the Hermitian Jacobi process, from a FFP theory perspective. 
For this matrix-valued process and the $\beta$-Jacobi extrapolation of its eigenvalue process, the following two asymptotic regimes are well-known:
\begin{enumerate}
\item For any fixed $\beta > 0$, the $\beta$-Jacobi process converges weakly as time $t \rightarrow +\infty$ to the so-called
$\beta$-\emph{Jacobi ensemble} whose density function is given by \cite{MR2719370}: 
$$
\prod_{i < j } |\lambda_i-\lambda_j|^{\beta} \prod_{i=1}^m \lambda_i^{p-1}(1-\lambda_i)^{q-1} d\lambda_i,
$$
where $p,q$ are positive parameters and $(\lambda_1,\dots, \lambda_m) \in [0,1]^m$. When $\beta \in \{1,2,4\}$, this density describes the joint distribution of the eigenvalues of a matrix drawn from the so-called Beta ensemble over the real, the complex and the quaternion algebras respectively. 
\item If time $t > 0$ is now fixed, then the empirical measure of the $2$-Jacobi process converges in the large-size limit (after the appropriate rescaling of the parameters) to the one-dimensional marginals of the free Jacobi process, a family of bounded operators initially defined and studied in \cite{MR2384475}. 
\end{enumerate}
The second convergence should hold for any fixed $\beta$, though the authors were not able to find it in the literature (see \cite{TT} for the case where $\beta$ varies and is comparable with the dimension). Later on, we will see that (Theorem \ref{eqn:asymptoticregime}) this result also holds when we further let $\beta = \infty$. 
Besides, we can interchange those two limits: the weak limit of the $2$-Jacobi ensemble is the free beta distribution (also called Wachter distribution \cite{Wac}). This is also the long-time limit of the one-dimensional marginal of the free Jacobi process or, equivalently, the spectral distribution of the angle operator of two free orthogonal projections (\cite{VDN}). It has the following explicit expression (using the name convention for the parameters given further below): 

\begin{equation*}
\nu_{\infty}^{\lambda,\theta}
= \frac{\bigl(\lambda-1\bigr)_{+}}{\lambda}\,\delta_{0}
+ \frac{\bigl((1+\lambda)\theta-1\bigr)_{+}}{\lambda\theta}\,\delta_{1}
+ \frac{\sqrt{(r_{+}-x)(x-r_{-})}}{2\pi\lambda\theta\,x(1-x)}\,
  \mathbf{1}_{[r_{-},r_{+}]}(x)\,dx,
\end{equation*}

\noindent where
\begin{equation*}
r_{\pm}
= (1+\lambda)\theta - 2\lambda\theta^{2}
\;\pm\; 2\theta\sqrt{\lambda\,(1-\theta)(1-\lambda\theta)}\,,
\end{equation*}

More generally, the moments of the Hermitian Jacobi process (so $\beta=2$) as well as their large-size limits were computed in the series of papers \cite{MR3842169, MR4173022,MR4954166}. Moreover, series expansions of $R$ and the free $S$-transforms of the free Jacobi process were derived in \cite{MR3679612} using the Lagrange inversion formula.

\subsection{Reminder on Jacobi processes and main results}
Let us state the defining equations of the $2$-Jacobi process and present our first main result. 
Let $(Y_t(d))_{t \geq 0}$ be a unitary Brownian motion of size $d$ and variance $t$ \cite{MR1430721}. Let
$P_m$ and $Q_p$ be two orthogonal projections of rank $m$ and $p$ respectively,  with $m \leq p$. The \emph{Hermitian Jacobi process} $(X_t^{(r,s,m)}, {t \geq 0})$ with parameters $(r,s,m)$ is the $m\times m$ matrix-valued stochastic process defined by:
\begin{equation}
\label{eqn:hermitianjacobi}
	X^{(r,s,m)}_t \oplus 0_{d-m} = P_mY_t(d)Q_p(Y_t(d))^\star P_m
\end{equation}
where
\begin{equation}
\label{eqn:parametersrsm}
	r:= p-m, \quad s := d-p-m:= q-m.
\end{equation}

\newcommand{\lambdapri}{\lambda_t^{(i,r,s,m)}}
\newcommand{\lambdaprj}{\lambda_t^{(j,r,s,m)}}

The eigenvalues process $(\lambda_t^{(i,r,s,m)}, t \geq 0)_{1 \leq i \leq m}$ of $(X^{(r,s,m)}_t, {t \geq 0})$ solves the following stochastic differential system:
\begin{align}\label{SDE}
	d\lambdapri &= \sqrt{2\lambdapri(1-\lambdapri)}dB_t^i \\ 
    &\hspace{0.75cm}+ \left[
		p-(p+q)\lambdapri + \sum_{j \neq i}\frac{\lambdapri(1-\lambdaprj) + \lambdaprj(1-\lambdapri)}{\lambdapri-\lambdaprj}\right] dt \nonumber,
\end{align}
where $(B_t^i, t \geq 0)_{1 \leq i \leq m}$ are independent Brownian motions \cite{MR3842169}. Up to the time change $t \mapsto 2t$, this random particle system is an instance of the $\beta$-Jacobi process studied in \cite{MR2384475} corresponding to the value $\beta = 2$. In particular, Corollary 2.1 there shows that the eigenvalues process is globally defined in time provided that $p \wedge q > m-1/2$. Corollary 7 in \cite{MR3296535} shows that the eigenvalues of $X^{(r,s,m)}_t$ are almost surely distinct and that global existence holds even under the weaker condition $p \wedge q > m-1$. Consequently, Theorem 4.1 in \cite{MR2384475} extends to $p \wedge q > m-1$ so that the probability distribution of the eigenvalues process $(\lambda_t^i, t \geq 0)_{1 \leq i \leq m}$ has a density with respect to Lebesgue measure in $\mathbb{R}^m$, for any $p \wedge q > m-1$. 

The $\beta$-Jacobi process $(\lambda_t^{(i,\beta,r,s,m)},t\geq 0)_{1\leq i\leq m}$ is the strong solution to the following stochastic differential equation : 
\newcommand{\lambdaprbetai}{\lambda_t^{(i,\beta,r,s,m)}}
\newcommand{\lambdaprbetaj}{\lambda_t^{(j,\beta, r,s,m)}}
\begin{align}\label{betaSDE}
	d\lambdaprbetai &= \sqrt{2\lambdaprbetai(1-\lambdaprbetai)}dB_t^i \\ 
    &\hspace{-1.25cm}+ \frac{\beta}{2}\left[
		p-(p+q)\lambdaprbetai + \sum_{j \neq i}\frac{\lambdaprbetai(1-\lambdaprbetaj) + \lambdaprbetaj(1-\lambdaprbetai)}{\lambdaprbetai-\lambdaprbetaj}\right] dt \nonumber.
\end{align}
We refer the reader to \cite{MR2719370} for conditions on $r,s,m$ and $\beta$ ensuring the existence of a strong solution to \eqref{betaSDE}.
\begin{remark}
\label{remark:marcus}
    Let $M_1 \in \mathcal{M}_{n_1,p}$ and $M_2 \in \mathcal{M}_{n_2,p}$ be $n_1 \times p$ and $n_2 \times p$ Gaussian matrices respectively. 
    The eigenvalues of 
\begin{equation*}
(M_1M^{\star}_1+M_2M^{\star}_2)^{-1/2}M_1M^{\star}_1(M_1M^{\star}_1+M_2M^{\star}_2)^{-1/2}
\end{equation*}
are distributed according to a $2$-Jacobi ensemble (matrix-variate beta law, \cite{collins2004}). In \cite{FFPP}, section 5.1., 
the author replaces the Gaussian matrices by Brownian matrices and observes that the eigenvalue process does not satisfy an autonomous SDE.
\end{remark}
The averaged characteristic polynomial of the $\beta$-Jacobi process is:
\begin{equation}
\label{eqn:charatersticpoly}
	\chi_t^{(r,s,m)}(x) := \mathbb{E}\prod_{i=1}^m(x-\lambdaprbetai),
\end{equation}
\newcommand{\charpoly}{\chi_t^{(r,s,m)}}
We intentionally do not record in the notation the dependence on $\beta$, since $\chi_t^{(r,s,m)}$ does not depend on this parameter as shown in \cite[Corollary 3.4]{MR4092353}.
In the series of papers \cite{MR4092353}, \cite{MR4673381}, \cite{MR4781073}, the roots of the characteristic polynomial of the (rescaled) $\beta$-Jacobi particle system valued in $[-1,1]$ are studied, among others. 
In particular, Theorem 4.4 from \cite{MR4781073} describes the dynamics of the roots of the averaged characteristic polynomial through an inverse Jacobi-heat equation.
Though the $\beta$-Jacobi particle systems in $[0,1]$ and in $[-1,1]$ are related by an affine transformation, the relation between elementary symmetric polynomials in both corresponding roots is far from being trivial. Nonetheless, one expects that the dynamics of the roots associated with the Jacobi particle system in $[0,1]$ are still governed by an inverse Jacobi-heat equation; this is indeed the case as we show later in Proposition \ref{prop:heatequationjacobi}, and yields the following result, Theorem \ref{cor:expressionaveragedcharacteristic}. 

\subsection{Frozen Jacobi process}
In the following Theorem, we expand the polynomial $\charpoly$ in the basis $(Q_j^{(r,s)})_{0\leq j \leq m}$ of the Jacobi polynomials defined by:
\begin{equation*}
	Q_j^{(r,s)}(x) := \frac{(r+1)_j}{j!}{}_2F_1(-j, r+s+j+1, r+1; x),
\end{equation*}
where ${}_2F_1$ stands for the Gauss hypergeometric function \cite[Chapter 4]{duverney2024introduction}. 

\begin{theorem}
	\label{cor:expressionaveragedcharacteristic}
	With the parameters $r,s,m$ defined as in equation \eqref{eqn:parametersrsm} above, let
	\begin{equation*}
		\chi_0^{(r,s,m)}(x) = \sum_{j=0}^m c_j^{(r,s,m)} Q_j^{(r,s)}(x),
	\end{equation*}
	be the expansion of the averaged characteristic polynomial at $t=0$ in the Jacobi polynomial basis. Then, for any time $t > 0$ :
	\begin{equation*}
		\chi_t^{(r,s,m)}(x) = \sum_{j=0}^m c_j^{(r,s,m)} e^{-(m-j)(r+s+1+m+j)t}Q_j^{(r,s)}(x).   \end{equation*}
	In particular, if $\chi_0^{(r,s,m)}(x) = (x-1)^m$, then
	\begin{equation*}
		c_j^{(r,s,m)} =  \frac{(r+s+1+2j)\Gamma(r+s+1+j)}{\Gamma(r+s+2+m+j)} \binom{m+s}{m-j}.
	\end{equation*}
\end{theorem}

Our second main result is concerned with the asymptotics of the root distribution of $\charpoly$: for any $t > 0$, it converges toward the one-dimensional marginal at time $t$ of the free Jacobi Process. We will give two proofs of this result. One is more technical and uses the notion of a finite $T$-transform; however, it is worth deriving the result this way, since it requires a technical lemma about the finite free $T$-transform, which is interesting in its own right. Moreover, this derivation also works in the broader setting where the roots of the sequence of polynomials we consider may not be confined to a compact set. 

Before stating our second theorem, let us recall a few facts on the free Jacobi process.
The free Jacobi process was introduced in \cite{MR2384475} as the large-size limit of the Hermitian Jacobi process $(X^{(r,s,m)}_{t/d(m)}, t \geq 0)$ in the asymptotic regime:
\begin{align}
	\label{eqn:asymptoticregime}
	\lambda = \lim_{m \rightarrow +\infty} \frac{m}{p(m)} > 0, \quad \lambda\theta = \lim_{m \rightarrow +\infty} \frac{m}{d(m)} \in (0,1].
\end{align}
or else, by using the parameters $r,s,m$:
$$
\lim_{m \rightarrow \infty} \frac{r(m)}{m} = \frac{1-\lambda}{\lambda},\quad 
\lim_{m \rightarrow \infty} \frac{s(m)}{m} = \frac{1-\theta(1+\lambda)}{\lambda \theta},
$$
More precisely, for any tuple of times $(t_1,\ldots,t_p)$, one has:
\begin{equation}
\label{eqn:convergenceofmoments}
 \frac{1}{m}{\rm Tr}(X^{(r(m),s(m),m)}_{t_1}\cdots X^{(r(m),s(m),m)}_{t_p}) \to_{m,+\infty} \tau(X_{t_1}^{(\lambda,\theta)} \cdots X_{t_p}^{(\lambda,\theta)})   
\end{equation}
where, for any time $t > 0$, the process $X_t^{(\lambda,\theta)}$ is defined by :
\begin{equation*}
X_t^{(\lambda,\theta)}:= PY_tQY_t^{\star}P    
\end{equation*}
where $P,Q$ are free (in Voiculescu's sense) orthogonal projections in a non commutative probability space $(\mathcal{A}, \tau)$ with ranks $\tau(P) = \lambda \theta$ and $\tau(Q) = \theta$, and $(Y_t)_{t \geq 0}$ is the so-called free unitary Brownian motion \cite{MR1430721}. In particular, under the asymptotic regime \eqref{eqn:asymptoticregime}, the convergence \eqref{eqn:convergenceofmoments} yields 
\begin{equation}
\label{eqn:convegencebetaonedimensional}
\frac{1}{m}\sum_{k=1}^m (\lambda_t^{(i,r(m),s(m),m)})^k \to \nu^{(\lambda, \theta)}_t(x^k)
\end{equation}
where $\nu_t^{(\lambda,\theta)}$ is the law of the one-dimensional marginal of the free Jacobi process at time $t$ and  $\nu^{(\lambda, \theta)}_t(x^k)$ is its $k$-th moment.

We will use $(x^{(r,s,m)}_t, t \geq 0)_{1 \leq i \leq m}$ for the (deterministic) process in $\mathbb{R}^m$ of the roots of $\charpoly$ (in the increasing order). The following Theorem extends the convergence \eqref{eqn:convegencebetaonedimensional} to the case $\beta = +\infty$. 
We denote by $\mu_t^{(r,s,m)}$ the counting measure of $x^{(r,s,m)}_t$ :
\begin{equation*}
	\mu_t^{(r,s,m)}:=\frac{1}{m}\sum_{i=1}^m\delta_{x_i(t)} = \mu\llbracket \charpoly \rrbracket
\end{equation*}

\begin{theorem}
\label{thm:maintheorem}
Assume $\mu_0^{(r(m),s(m),m)}$ converges weakly in the regime specified in \eqref{eqn:asymptoticregime} to $\nu_0^{(\lambda,\theta)}$. Then, for any time $t > 0$, the counting measure $\mu_t^{(r,s,m)}$ of $x^{(r,s,m)}_t$ converges to $\nu^{(\lambda,\theta)}_t$, which is equivalent to:

$$
\frac{1}{m} \sum_{k=1}^m (x_t^{(r,s,m)})^k \to \nu_t^{(\lambda,\theta)}(x^k),\quad k\geq 1.
$$
\end{theorem}

A direct way of proving Theorem \ref{thm:maintheorem} is to find a differential system satisfied by the moments of $\mu^{(r,s,m)}_t$ and to take the large-size limit there. Then, we recognize the limiting ODE as the one satisfied by the moments of the one-dimensional marginal of the free Jacobi process. This line of argument works well because we already have a candidate probability distribution whose moments satisfy the limiting ODE, and since the involved distributions are compactly supported for any time $t$. 
In a non compact/unbounded setting, it would be far more difficult to use a similar argument. For that reason, we develop a different strategy making use of the $S$ or the $T$ transforms as explained in the next paragraph. 

\subsection{PDE for the $S$ and $T$ of the free Jacobi process}
The first definition of the finite free $S$-transform (which is in turn related to the finite free $T$-transform) is due to Marcus \cite{MR4408504}. Few years later, another definition was proposed
and was motivated by the CLT for finite free multiplicative convolution \cite{arizmendi2024s}. Marcus's definition is general enough to cover the finite free multiplicative convolutions of distributions supported on the positive half-line and on the circle, but is much less practical for computations. For this reason, we prefer to work with the definition proposed in \cite{arizmendi2024s}.

On the free probabilistic side, recall that if $\mu \neq \delta_0$ is a probability measure on $\mathbb{R}_+$ and if  
\begin{equation*}
\psi_{\mu} \colon \mathbb{C}\setminus \mathbb{R}_+ \to \mathbb{C}    
\end{equation*} 
is the moment generating function of $\mu$, then $\psi_{\mu}$ is univalent on the left half-plane $i\mathbb{C}_+$. Besides $\psi_{\mu}(i\mathbb{C}_+)$ is a domain contained in the disc with diameter $(-1+\mu\{0\}, 0)$ and contains indeed this interval. The free $S$-transform of $\mu$ by, for $z$ in the vicinity of $(-1+\mu\{0\}, 0)$ is then defined by
$$
\mathcal{S}_{\mu}(z) = \frac{1+z}{z} \psi_{\mu}^{(-1)}(z).
$$

Its main property is the factorization under the free multiplicative convolution of distributions on the positive half-line :
$$
\mathcal{S}_{\mu_1 \boxtimes \mu_2}(z) = \mathcal{S}_{\mu_1}(z)\mathcal{S}_{\mu_2}(z).
$$
\newcommand{\Sfree}{\mathcal{S}_t^{(\lambda,\theta)}}
\newcommand{\Tfree}{\mathcal{T}_t^{(\lambda,\theta)}}

Hereafter, we shall denote $\mathcal{S}_t^{(\lambda,\theta)}$ the $S$-transform of $\nu_t^{(\lambda,\theta)}$ at time $t$. By using free stochastic calculus, we will derive a partial differential equation satisfied by $\Sfree$. The free $T$-transform of $\mu$ is defined as the reciprocal of its $S$-transform : 
$$
\mathcal{T}_{\mu}(z) := \frac{1}{\mathcal{S}_\mu(z)},
$$
and we similary denote $\mathcal{T}_t^{(\lambda, \theta)}$ the $T$-transform of $\nu_t^{(\lambda,\theta)}$. 

\begin{theorem}
	\label{thm:derivationfreestransform}
	The $S$-transform of the free Jacobi process $\Sfree$ satisfies, locally around the origin, the PDE:
	\begin{align}\label{PDE-S}
		&\partial_t \mathcal{S}_t^{(\lambda,\theta)}(z) =  (2\lambda \theta z+1) \mathcal{S}_t^{(\lambda,\theta)}
		-\theta (1+2\lambda z)[\mathcal{S}_t^{(\lambda,\theta)}(z)]^2 -\frac{\theta z (1+\lambda z)}{2}\partial_z[\mathcal{S}_t^{(\lambda,\theta)}(z)]^2, 
        \\ 
        &\mathcal{S}_0^{(\lambda,\theta)}(z) = 1. \nonumber
	\end{align} 
    Equivalently,
    \begin{align}\label{PDE-T}
    	&\partial_t \mathcal{T}_t^{(\lambda,\theta)}(z) = (2(1-z)\lambda\theta -1)\mathcal{T}_t^{(\lambda,\theta)}(z) + \theta(1-2\lambda (1-z)) \\
		&\hspace{7cm}+ \theta(1-\lambda{(1-z)})(1-z){\partial_{z}\log \mathcal{T}_t^{(\lambda,\theta)}(z)},\nonumber \\
		&\mathcal{T}_0^{(\lambda,\theta)}(z) = 1. \nonumber
    \end{align}
\end{theorem}
\begin{remark} 
Since $\nu_t^{(\lambda,\theta)}$ is compactly supported on $[0,1]$, both its $S$ and its $T$-transforms are analytic around the origin and in a band comprising the segment $(-1+\mu(0),0)$. 
The above PDEs satisfied by $\mathcal{S}_t^{(\lambda,\theta)}$ and by $\mathcal{T}_t^{(\lambda,\theta)}$ yields triangular systems of ordinary differential equations on their Taylor coefficients, which can be shown to admit unique solutions. Hence, by standard arguments, we can argue that $\Sfree$ and $\Tfree$ are the unique solutions to \eqref{PDE-S} and \eqref{PDE-T} 
which are analytic around the origin.
\end{remark}
We now introduce the finite free $S$ and $T$-transforms and state the main technical result that allows the derivation of the PDEs \eqref {PDE-S} and \ref{PDE-T} without appealing to free stochastic calculus. Let $p$ be a monic polynomial with degree $m$:
$$
p_m(x) = \sum_{k=0}^m (-1)^k e^{(m)}_k(p) x^{m-k}, \quad e^{(m)}_0(p) = 1.
$$
Let $r$ the multiplicity of $0$ in $p_m$ and assume that $p_m$ \emph{has only positive real roots}.
Then the \emph{finite free $S$ transform of} $p_m$, denoted hereafter by $S^{(m)}_{p_m},$ is the function defined on the set of points $\{-k/m, k \in \{1,m-r\} \}$ defined by \cite{arizmendi2024s}:
$$
{S}_{p_m}^{(m)}\left(-\frac{k}{m}\right) := \frac{m-k+1}{k}\frac{e^{(m)}_{k-1}(p)}{e^{(m)}_k(p)},\quad k \in \{1,m-r\}.
$$
We find it also convenient to introduce the finite free $T$ transform of $p_m$ : it is the piecewise right-continuous function defined on the open interval $(0,1)$ by: 
$$
T_{p_m}^{(m)}(v) = \begin{cases}
                    0, & v \in (0,\frac{r}{m})\\
\displaystyle                    \frac{m-k+1}{k}\frac{e^{(m)}_{m-k+1}(p)}{e^{(m)}_{m-k}(p)}, & v \in [\frac{k-1}{m}, \frac{k}{m}), \quad k = r+1, \ldots, m.
                 \end{cases}
$$ 
Equivalently, 
\begin{equation}\label{Defi-T-Trans} 
T_{p_m}^{(m)}(v) = \begin{cases}
                    0, & v \in (0,\frac{r}{m})\\
\displaystyle                    \frac{m-\ceil{mv}+1}{\ceil{mv}}\frac{e^{(m)}_{m-\ceil{mv}+1}(p)}{e^{(m)}_{m-\ceil{mv}}(p)}, & v \in [(r/m),1),
                 \end{cases}
\end{equation}
where $\ceil{\cdot}$ is the ceiling function. 
Note in passing that while the finite free $S$-transform is defined on a finite set of real numbers, the finite free $T$-transform is defined in the open interval $(0,1)$.  Nonetheless, for any $k \in \{1, \dots, m-r\},$ 
\begin{equation}
\label{eqn:relation}
   T_{p_m}^{(m)}\left(\frac{m-k}{m}\right) = \frac{1}{\displaystyle {S}_{p_m}^{(m)}\left(-\frac{k}{m}\right)}.
\end{equation}

Let $\nabla^{(m)}$ be the operator of finite right-differentiation with step $1/m$ acting on functions of $\mathbb{R}$:
$$ 
	\nabla^{(m)}\colon \mathbb{R}^{\mathbb{R}}\to \mathbb{R}^{\mathbb{R}},\quad \nabla^{(m)} g (v) := m\left(g \left(v+\frac{1}{m}\right) - g (v)\right),\quad v \in \mathbb{R}
.$$

Likewise, the discrete left-derivative of $g$ is 
$$
\tilde{\nabla}^{(m)}\colon \mathbb{R}^{\mathbb{R}} \to \mathbb{R}^{\mathbb{R}}, \quad \tilde{\nabla}^{(m)}g(v) = m\left(g(v)-g\left(v-\frac{1}{m}\right)\right),\quad v \in \mathbb{R}.
$$
Most of the following statements will be formulated by using the right-derivative operator $\nabla^{(m)}$, although all have their counterparts valid for the left derivative operator $\tilde{\nabla}$.

Before stating the next Theorem \ref{thm:finitefreedifference}, recall that if $(p_m)_{m\geq 1}$ is a sequence of polynomials with increasing degree $m$, we denote by $\mu\llbracket p_m \rrbracket $ the counting measure of its roots : 
$$
\mu\llbracket p_m\rrbracket : = \frac{1}{m}\sum_{z\,:\,P(z)=0} \delta_z.
$$
If $\mu\llbracket p_m \rrbracket$ converges weakly to some probability measure $\mu$, we say that $(p_m)_{m \geq 1}$ is a \emph{converging sequence}. 
From \cite[Theorem 1.1]{arizmendi2024s}, it follows that if for any $m \geq 1$, $p_m$ has only positive roots and if $\mu \neq \delta_0$, 
then $\mu\llbracket p_m \rrbracket $ converges weakly to $\mu$ if and only if 
$(S^{(m)}_{p_m})_m$ converges to the free $S$-transform $\mathcal{S}_\mu$ of $\mu$ in the following sense: for any $v \in (0,1-\mu\{0\})$ 
and any sequence $(k_m)_{m \geq 1}$ such that
\begin{equation*}
 1 \leq k_m \leq m, \quad \lim_{m \rightarrow +\infty} \frac{k_m}{m} = v,    
\end{equation*}
one has: 
\begin{equation*}
\lim_{m \rightarrow +\infty} {S}_{p_m}^{(m)}\left(-\frac{k_m}{m}\right) = S_{\mu}(-v).     
\end{equation*}
We will also make use of the following interpolation over the interval $(0,1)$ of the finite free $S$-transform of the polynomials $p_m$:
	$$
		\mathcal{S}^{(m)}_{p_m} \colon (0,1)\to \mathbb{R},\quad \mathcal{S}^{(m)}_{p_m}(v) = S^{(m)}_{p_m}(-\frac{\ceil{mv}}{m}),\quad v \in \left(0,1-\frac{r_m}{m}\right),
	$$
	where we recall that $\ceil{\cdot}$ is the ceiling function and $r_m$ is the multiplicity of $0$ in $p_m$ (not to be confused with $r(m) = p(m) - m)$. In this regard, we again
emphasize that this interpolated $S$ transform $\mathcal{S}^{(m)}_{p_m}$ is only defined on the interval $(0,1-r_m/m)$, while $T^{(m)}_{p_m}$ is defined on the whole interval $(0,1)$. 
The relation \eqref{eqn:relation} then becomes 
$$
T^{(m)}_{p_m}(1-v) = \frac{1}{\mathcal{S}^{(m)}_{p_m}(v)}, \quad v \in\left(0, 1-\frac{r_m}{m}\right).
$$

\begin{theorem}
	\label{thm:finitefreedifference}
	Let $(p_m)_{m\geq 1}$ be a converging sequence of monic polynomials with positive roots and increasing degrees $m$. Let $\mu$ be the limiting measure :
	$$
		\mu\llbracket p_m \rrbracket \underset{m,+\infty}{\longrightarrow} \mu,
	$$
	and assume that $\mu \neq \delta_0$. Then, 
	\begin{align*}
		[\nabla^{(m)}](\mathcal{S}^{(m)}_{p_m})(v) = -\partial_v{\mathcal{S}}_{\mu}(-v) + o(1), \quad m \rightarrow +\infty,
	\end{align*}
	locally uniformly on $v \in (0,1-\mu(\{0\})$.
\end{theorem}
\begin{remark}
 From the statement of Theorem \ref{thm:finitefreedifference}, it is expected that the iterated discrete derivatives of the finite $S$-transform converge to the iterated derivatives of the $S$- transform. We will not prove this assertion; a possible strategy would be to use the same arguments as exposed in the proof of Theorem \ref{thm:finitefreedifference} in a recursive way. 
\end{remark}

\begin{corollary}
Continuing the setting of Theorem \ref{thm:finitefreedifference}, then for any $v \in (\mu{(\{ 0 \})},1)$, 
$$
[\nabla^{(m)}](T^{(m)}_{p_m})(v) = \partial_v \mathcal{T}_\mu(v-1) + o(1)
$$
as $m$ goes to infinity.
\end{corollary}
Hereafter, we denote by $\mathcal{S}^{(r,s,m)}_t$ (resp. $T^{(r,s,m)}_t$) the finite free $S$ (resp. the finite free $T$) transform of the counting measure $\mu^{(r,s,m)}_t$ of the roots 
$x^{(r,s,m)}_t$ of $\chi^{(r,s,m)}_t$. 
In Proposition \ref{cor:odefiniteTtransform}, we exhibit difference-differential equations satisfied by the finite free $T$ transforms of $\mu^{(r,s,m)}_t$ and pass to the limit as $m$ goes to infinity (using Theorem \ref{thm:finitefreedifference}) to recover the PDEs stated in Theorem \ref{thm:derivationfreestransform}.

In addition to this Section, the paper has three other Sections. In Section \ref{sec:stransformfree}, we prove Theorem \ref{thm:derivationfreestransform}, we provide a formula for the stationary analytic solution of the PDE satisfied by the $S$-transform in Section \ref{ssec:stationnarysolution},  in Section \ref{ssec:charcurves}, we analyse the characteristic curves. 

In Section \ref{sec:averaged}, we prove Theorem \ref{cor:expressionaveragedcharacteristic}. In Section \ref{ssec:unitaryhermite}, the relation between the characteristic polynomial of the Hermitian Jacobi process and the Hermite unitary polynomial (the characteristic polynomial of the Unitary Brownian motion) is elucidated. In Section \ref{ssec:anotherformula}, we derive the formula for the characteristic polynomial of the Jacobi process without relying on the dynamics of its roots. 

In Section \ref{sec:frozen}, we prove Theorem \ref{thm:maintheorem} by computing the asymptotics of the moments of the frozen Jacobi process. In Section \ref{sec:finitefreediff} we prove the technical statement, Theorem \ref{thm:finitefreedifference}, about the convergence of the finite free differences of the finite free $S$ transform and use it to rederive, in Theorem \ref{thm:PDEfromfinitefree}, the PDE satisfied by the free $S$ transform of the Jacobi process

\section{\texorpdfstring{$S$}{S} transform of the free Jacobi Process}
\label{sec:stransformfree}
This section is devoted to the study of the $S$-transform of the free Jacobi process and to proving Theorem \ref{thm:derivationfreestransform}. Once proved, it will be shown that the unique non trivial stationary solution are coincides with the $S$ transform of the angle operator $PQP$ between two free orthogonal projections (viewed as an operator in the compressed probability space). Though this coincidence is expected because the free Jacobi process converges weakly to $PQP$, it is not straightforward from a partial differential equation perspective. 

To proceed, recall from \cite{MR2384475} that the Cauchy transform $G_t^{(\lambda,\theta)}$ of the spectral distribution $\nu_t^{(\lambda,\theta)}$ of the free Jacobi process as time $t$ satisfies the following non-linear partial differential equation (PDE):
\begin{equation*}
	\partial_tG_t^{(\lambda,\theta)} = \partial_z \left\{
	[(1-2\lambda \theta) z -\theta(1-\lambda)]G_t^{(\lambda,\theta)}
	+\lambda \theta z(z-1)[G_t^{(\lambda,\theta)}]^2\right\}.
\end{equation*}
Consequently, the moment generating function
\begin{equation*}
	\mathcal{M}_t^{(\lambda,\theta)}(z) = \frac{1}{z}G_t^{(\lambda,\theta)}\left(\frac{1}{z}\right)
\end{equation*}
satisfies:
\begin{equation}
\label{PDE-MGF}
	\partial_t\mathcal{M}_t^{(\lambda,\theta)} = -z \partial_z \left\{
	[(1-2\lambda \theta)  -\theta(1-\lambda)z]\mathcal{M}_t^{(\lambda,\theta)}
	+\lambda \theta (1-z)[\mathcal{M}_t^{(\lambda,\theta)}]^2\right\},
\end{equation}
which reduces when $\lambda =1, \theta =1/2,$ to 
\begin{equation*}
	\partial_t\mathcal{M}_t^{(1,1/2)} = -\frac{z}{2} \partial_z \left\{(1-z)[\mathcal{M}_t^{(1,1/2)}]^2\right\}.
\end{equation*}
The analytic (around the origin) solution of this PDE subject to the initial data 
\begin{equation*}
\mathcal{M}_0^{(1,1/2)}(z) = \frac{1}{1-z}     
\end{equation*}
was determined in \cite{MR3071702} and used in \cite{MR3679612} to derive expansions of the corresponding $R$ and $S$ transforms, based on Lagrange inversion formula. Other degenerate limiting regimes were investigated in \cite{MR4751304} for the Jacobi eigenvalue process, and gave rise to Marchenko-Pastur and Wigner distributions. 
At the root level, these regimes correspond to the ones studied in \cite{MR1342385}.

\subsection{Proof of Theorem \ref{thm:derivationfreestransform}}
For any time $t \geq 0$, denote $\eta_t^{(\lambda,\theta)}$ the local inverse around zero of 
\begin{equation*}
\psi_t^{(\lambda,\theta)} := \mathcal{M}_t^{(\lambda,\theta)}-1,
\end{equation*}
Since $\psi_t^{(\lambda,\theta)}[\eta_t^{(\lambda,\theta)}(z)] = z$, then the chain rule entails:
	\begin{eqnarray*}
		1 &=& \partial_z\psi_t^{(\lambda,\theta)}[\eta_t^{(\lambda,\theta)}(z)]\partial_z\eta_t^{(\lambda,\theta)}(z)
		\\
		0 &=& \partial_t\psi_t^{(\lambda,\theta)}[\eta_t^{(\lambda,\theta)}(z)] + \partial_t[\eta_t^{(\lambda,\theta)}](z)
		\partial_z\psi_t^{(\lambda,\theta)}[\eta_t^{(\lambda,\theta)}(z)].
	\end{eqnarray*}
	It follows that
	\begin{align*}
		 & \partial_t[\eta_t^{(\lambda,\theta)}]  =  -\partial_t\psi_t^{(\lambda,\theta)}[\eta_t^{(\lambda,\theta)}] \partial_z[\eta_t^{(\lambda,\theta)}]
		\\& =  [\eta_t^{(\lambda,\theta)}]  \partial_z \left\{
		[(1-2\lambda \theta)  -\theta(1-\lambda)z](\psi_t^{(\lambda,\theta)}+1)
		+\lambda \theta (1-z)[(\psi_t^{(\lambda,\theta)}+1)]^2\right\}
		[\eta_t^{(\lambda,\theta)}]                                                                                                                        \\& \times \partial_z[\eta_t^{(\lambda,\theta)}]
		\\& = [\eta_t^{(\lambda,\theta)}]
		\partial_z \left\{(1+z)\left[(1-2\lambda \theta)  -\theta(1-\lambda)\eta_t^{(\lambda,\theta)}\right]
		+\lambda \theta (z+1)^2 \left(1-\eta_t^{(\lambda,\theta)}\right)\right\}
		\\& = [\eta_t^{(\lambda,\theta)}]
		\partial_z \left\{(1-\lambda \theta +z +\lambda \theta z^2)
		-\theta(1+(1+\lambda)z +\lambda z^2)\eta_t^{(\lambda,\theta)}\right\}
		\\& = (2\lambda \theta z+1)[\eta_t^{(\lambda,\theta)}]
		-\theta(1+\lambda +2\lambda z) \left[\eta_t^{(\lambda,\theta)}\right]^2
		-\frac{\theta}{2}(1+z)(1+\lambda z) \partial_z \left[\eta_t^{(\lambda,\theta)}\right]^2.
	\end{align*}
	Multiplying both sides by $(1+z)/z$ and using the fact that
	\begin{align*}
		\partial_z\left[\mathcal{S}_t^{(\lambda,\theta)}\right]^2 & = \frac{(1+z)^2}{z^2}\partial_z \left(\eta_t^{(\lambda,\theta)}\right)^2
		-2\frac{1+z}{z^3}\left(\eta_t^{(\lambda,\theta)}\right)^2
		\\& = \frac{(1+z)^2}{z^2}\partial_z \left(\eta_t^{(\lambda,\theta)}\right)^2
		-\frac{2}{z(1+z)}\left(\mathcal{S}_t^{(\lambda,\theta)}\right)^2,
	\end{align*}
	we get the desired result.

\subsection{Stationary analytic (around the origin) solutions of \texorpdfstring{\eqref{PDE-S}}{}} 
\label{ssec:stationnarysolution}
These are the solutions of \eqref{PDE-S} which do not depend on time. As such, they satisfy the ordinary differential equation (hereafter ODE):    
\begin{equation*}
(2\lambda \theta z+1) \mathcal{S}^{(\lambda,\theta)}
		-\theta (1+2\lambda z)[\mathcal{S}^{(\lambda,\theta)}(z)]^2 -\frac{\theta z (1+\lambda z)}{2}\partial_z[\mathcal{S}^{(\lambda,\theta)}]^2(z)=0.
\end{equation*}
The zero function is obviously the trivial solution. Otherwise, any non-zero stationary analytic solution satisfies
\begin{equation*}
\theta z (1+\lambda z)\partial_z\mathcal{S}^{(\lambda,\theta)}(z) + \theta (1+2\lambda z)\mathcal{S}^{(\lambda,\theta)}(z) = 1+ 2\lambda \theta z.
\end{equation*}
Noting that the LHS of this ODE may be written as $\partial_z[\theta z (1+\lambda z)\mathcal{S}^{(\lambda,\theta)}(z)]$, we readily see that 
\begin{equation*}
  \mathcal{S}^{(\lambda,\theta)}(z)=  \frac {\lambda\theta z^2 + z + \theta c} {\theta z (\lambda z + 1 )},\quad z\in \mathbb{C}\setminus\{0,-\frac{1}{\lambda}\}\ {\rm and}\ c\in\mathbb{C}. 
\end{equation*}
But analyticity around the origin forces $c = 0$ and we end up (after removing the singularity at $z=0$) with 
\begin{equation*}
  \mathcal{S}^{(\lambda,\theta)}(z)=  \frac {\lambda\theta z^2 +  z } {\theta z (\lambda z + 1 )} = \frac {\lambda\theta z +  1 } { \lambda\theta  z + \theta}.
\end{equation*}
This is the $S$ transform $S_Q$ of an orthogonal projection $Q$ of rank $\theta$ taken at $\lambda \theta z$.  Equivalently, it is the $S$ transform of the weak limit as $t \rightarrow +\infty$ of the spectral
distribution of the free Jacobi process.
The latter may be defined as the spectral distribution of the angle operator $PQP$ between two free orthogonal projections $P$ and $Q$ in a non commutative probability space $(\mathscr{A}, \tau)$ with ranks $\tau(P)= \lambda \theta$ and $\tau(Q) = \theta$, viewed as an operator in the compressed probability space $(P\mathscr{A}P, \tau/\tau(P))$. 

To see this, we firstly use the Nica-Speicher convolution semi-group (see \cite{MR2266879}, Lecture 14) and deduce that the Voiculescu's $R$ transform $\mathcal{R}_{PQP}$ of $PQP$ in the compressed space is given by 
\begin{equation*}
R_Q(\lambda \theta z)     
\end{equation*}
where $R_Q$ is the Voiculescu's $R$ transform of $Q$ in $(\mathscr{A}, \tau)$, (see \cite{MR2266879}, page 211). Secondly, let 
\begin{equation*}
\mathcal{C}_{PQP}(z) := z \mathcal{R}_{PQP}(z), \quad C_Q(z) := zR_Q(z), 
\end{equation*}
be the shifted Voiculescu R transforms, also known as free cumulant generating functions (see \cite{MR2266879}, Lecture 16). Then, we equivalently get
\begin{equation*}
 \mathcal{C}_{PQP}(z) = zR_Q(\lambda \theta z) = \frac{1}{\lambda \theta}C_Q(\lambda \theta z). 
\end{equation*}
Consequently, the local inverse around $z=0$ of $\mathcal{C}_{PQP}$ is given by 
\begin{equation*}
\mathcal{C}_{PQP}^{-1}(z) = \frac{1}{\lambda \theta}C_Q^{-1}(\lambda \theta z).  
\end{equation*}
Finally, the functional relations (see \cite{MR2266879}, Lecture 18):
$$\mathcal{C}_{PQP}^{-1}(z) = z\mathcal{S}_{PQP}(z), \quad C_Q^{-1}(z) = zS_Q(z),
$$
yield 
\begin{equation*}
 z\mathcal{S}_{PQP}(z) = \frac{\lambda \theta z}{\lambda \theta}S_Q(\lambda \theta z) \quad \Leftrightarrow \quad \mathcal{S}_{PQP}(z) = S_Q(\lambda \theta z) = \mathcal{S}^{(\lambda,\theta)}(z).     
\end{equation*}
In a nutshell, 
\begin{equation*}
 \mathcal{S}^{(\lambda,\theta)}(z) = \lim_{t \rightarrow +\infty} \mathcal{S}_t^{(\lambda,\theta)}(z)   
:= \mathcal{S}_{\infty}^{(\lambda,\theta)}(z)
\end{equation*}
is the non-trivial analytic stationary solution of the PDE \eqref{PDE-S}.

\subsection{Analysis of the Characteristic curves for \texorpdfstring{$\lambda = 1, \theta=1/2$}{}}
\label{ssec:charcurves}
In this paragraph, we restrict our attention to the special parameter set \(\theta = 1/2\), \(\lambda = 1\), and assume for sake of simplicity that $\mathcal{S}^{(1,1/2)}_0(z) =1$ (in which case the initial spectral distribution of the free Jacobi process is $\delta_1$). This restriction is mainly motivated by the fact that the moment-generating function $\mathcal{M}_t^{(1,1/2)}$
admits an explicit expression, as proved in \cite{MR3071702}, from which a series expansion of $\mathcal{S}_t^{(1,1/2)}$ was derived in \cite{MR3679612} making use of Lagrange inversion formula. 
As we shall see below, the analysis of the Characteristic curves of the pde \eqref{PDE-S} in the special case corresponding to $\lambda = 1, \theta=1/2$ leads to a considerably more tractable expression of 
$\mathcal{S}^{(1,1/2)}_t$ than the one displayed in Proposition 5.2. from \cite{MR3679612}. 
\begin{proposition}
	\label{prop:characteristics}
	Let $t\geq 0$. The S transform $\mathcal{S}_t^{(1,1/2)}$ is given locally around $z=0$ by
	\begin{align}
        \label{eqn:stransformunundemi}
		\mathcal{S}_t^{(1,1/2)}(z) & = \frac{z^2 + 2z - \kappa_t^{-1}(z)}{z(1+z)} = \mathcal{S}_{\infty}^{(1,1/2)}(z) - \frac{\kappa_t^{-1}(z)}{z(1+z)},
	\end{align}
	where $\kappa_t$ is defined by
	$$
		\kappa_t(z) = \frac{(1+\sqrt{1+z})\xi_{2t}(\sqrt{1+z}) + (\sqrt{1+z}-1)}{1-\xi_{2t}(\sqrt{1+z})}.
	$$
	For any time $t > 0$, $\xi_{2t}$ is the inverse, in a vicinity of $u=1$, of the Herglotz transform of the spectral distribution of the free unitary Brownian motion at time $2t$ \cite{MR1430721}:
	$$
		\xi_{2t}(u)=\frac{u-1}{u+1} e^{ut}.
	$$
\end{proposition}
\begin{proof}
	If $\theta = 1/2, \lambda =1,$ then the PDE \eqref{PDE-S} reduces to:
	\begin{align}\label{PDE-S00}
		\partial_t \mathcal{S}_t^{(1,1/2)}(z) & =  (1+z)\mathcal{S}_t^{(1,1/2)}(z) -\frac{1+2z}{2}\left[\mathcal{S}_t^{(1,1/2)}(z)\right]^2
		-\frac{z(1+z)}{4}\partial_z\left[\mathcal{S}_t^{(1,1/2)}\right]^2(z).
	\end{align}
	For a fixed $z$ near the origin, the characteristic curve $s \mapsto z(s)$ of the PDE \eqref{PDE-S00} starting at $z$ satisfies locally (in time) the ODE
	\begin{equation}\label{charasteristic}
		z'(s)=\frac{1}{2}z(s)(1+z(s))f(s), \quad z(0)=z.
	\end{equation}
	where we set:  
    $$
    f(s) :=\mathcal{S}_t^{(1,1/2)}(z(s)).
    $$
	It follows that 
    \begin{equation}\label{characteristic1}
		f'(s)=(1+z(s))f(t)-\frac{1+2z(s)}{2}f(s)^2, \quad f(0)=\mathcal{S}_0^{(1,1/2)}(z).
	\end{equation}
	Substituting the ODE \eqref{charasteristic} into \eqref{characteristic1} yields
	\begin{equation*}
		2\frac{z''(s)z(s)(1+z(s))-z'(s)^2(1+2z(s))}{z(s)^2(1+z(s))^2}=\frac{2z'(s)}{z(s)}-2\frac{z'(s)^2(1+2z(s))}{z(s)^2(1+z(s))^2},
	\end{equation*}
	which reduces to
	\begin{equation*}
		z''(s)= z'(s)(1+z(s)).
	\end{equation*}
	A first integration gives
	\begin{equation*}
		z'(s) - z'(0) = z'(s) - \frac{z(1+z)}{2}f(0) = z(s) +\frac{z^2(s)}{2} - z - \frac{z^2}{2},
	\end{equation*}
	or equivalently
	\begin{equation*}
		z'(s)  = z(s) +\frac{z^2(s)}{2} - z\left(1-\frac{S_0^{(1,1/2)}(z)}{2}\right) - \frac{z^2}{2}(1-S_0^{(1,1/2)}(z)).
	\end{equation*}
	Specializing the last ODE to $\mathcal{S}_0^{(1,1/2)}(z) = 1$, a second integration then gives (locally in time)
	\begin{equation*}
		\frac{z(s)+1-\sqrt{1+z}}{z(s) + 1+ \sqrt{1+z}} = \frac{\sqrt{z+1}-1}{\sqrt{z + 1} + 1} e^{\sqrt{1+z} s} = \xi_{2s}(\sqrt{1+z}),
	\end{equation*}
	while equation \eqref{charasteristic} yields again:
	\begin{equation*}
		\mathcal{S}_t^{(1,1/2)}(z(s)) = f(s) = \frac{z^2(s) + 2z(s) - z}{z(s)(1+z(s))}.
	\end{equation*}
	Consequently, it holds locally in time that: 
	\begin{equation}\label{KappaMap}
		z(s) = \frac{(1+\sqrt{1+z})\xi_{2s}(\sqrt{1+z}) + (\sqrt{1+z}-1)}{1-\xi_{2s}(\sqrt{1+z})} = \kappa_s(z).
	\end{equation}
	
    Now, recall from \cite{MR1430721}, Lemma 12, that for any fixed time $t$, there exists a Jordan $\Gamma_t$ lying in the right half-plane and containing $w=1$ where $\xi_{t}$ is a one-to-one map from $\Gamma_t$ onto the open unit disc. As a matter of fact, there exists a neighbourhood of the origin such that for any $z$ there, the associated characteristic curve $s \mapsto z(s)$ is globally defined in time. 
    
    Finally, fix $t > 0$ and assume for a while that $z \mapsto \kappa_t(z)$ is locally invertible around the origin. Then $z = \kappa_t^{-1}(z(t))$ and we finally get:
	\begin{align*}
		\mathcal{S}_t^{(1,1/2)}(z(t)) & = \frac{z^2(t) + 2z(t) - \kappa_t^{-1}(z(t))}{z(t)(1+z(t))},
	\end{align*}
    whence the expression of $\mathcal{S}_t^{(1,1/2)}(z)$ follows. 
    
	Coming back to the local invertibility of $\kappa_t$ around the origin, observe that
	\begin{equation*}
		\kappa_t(z) = \sqrt{1+z} \frac{1+\xi_{2t}(\sqrt{1+z})}{1-\xi_{2t}(\sqrt{1+z})} - 1,
	\end{equation*}
	so that one only needs to check the same property for the map
	\begin{equation*}
		v \mapsto v\frac{1+\xi_{2t}(v)}{1-\xi_{2t}(v)},
	\end{equation*}
	around $v=1$. Since the derivative of this map at $v=1$ is given by $1+2\xi_{2t}'(1) = 1+e^t \neq 0$, the proposition is then proved.
\end{proof}
\begin{remark}
	For any time $t\geq 0$, let $V_t$ be defined by:
	\begin{equation*}
		V_t(z) = z(1+z) \mathcal{S}_t^{(1,1/2)}(z),
	\end{equation*}
	From the PDE \eqref{PDE-S00}, we readily infer
	\begin{equation*}
		\partial_t V_t = (1+z)V_t - \frac{1}{4}\partial_z V_t^2,
	\end{equation*}
	and in turn $z \mapsto \kappa_t^{-1}(z) = z^2 +2z -V_t(z)$ satisfies
	\begin{equation*}
		\partial_t\kappa_t^{-1} = \frac{1}{4}\partial_z [\kappa_t^{-1}]^2
		-\frac{z^2+2z}{2}\partial_z \kappa_t^{-1}.
	\end{equation*}
On the other hand, the proposition yields the limit
\begin{equation*}
 \lim_{t \rightarrow + \infty} \kappa_t^{-1} = 0,    
\end{equation*}
This aligns with the fact that the Jordan domain $\Gamma_t$ shrinks to $w=1$ in the limit as $t \rightarrow +\infty$ and the fact that $= \xi_{2t}(1) = 0 = \kappa_t(0)$ for any time $t \geq 0$.
\end{remark}

\section{Averaged characteristic polynomial of the Hermitian Jacobi Process}
\label{sec:averaged}
In this section, we proceed to the study of the averaged characteristic polynomial of the Hermitian Jacobi process; we prove, in particular, Theorem \ref{thm:maintheorem}. In Proposition \eqref{prop:heatequationjacobi}, we revisit Voit's result showing that the averaged characteristic polynomial satisfies the inverse Jacobi-heat equation for $p \wedge q > m-(1/2)$. 
In doing so, we prove it directly using stochastic calculus, without relying on elementary symmetric polynomials. We also emphasise that this result holds for the larger set $p \wedge q > m-1$ by virtue of Corollary 7 in \cite{MR3296535}. Proposition \ref{prop:positivityoftheroot}
shows that all its roots remain in the interval $(0,1)$. In this respect, we would like to stress that this last result is far from obvious from a differential-equation perspective, since the dynamics of the roots are highly nonlinear. Actually, it holds true in our framework since this deterministic particle system results from 'freezing' the eigenvalues process \eqref{SDE}.  Proposition \ref{prop:relationunitaryjacobi} is the finite free analogue of Corollary 3.3. proved in \cite{MR3071702}: when $\lambda =1, \theta=1/2$, the image of the roots of
the averaged characteristic polynomial of the Hermitian Jacobi process by the transformation $x \mapsto 2\arccos(\sqrt{x})$ coincides with the roots of the Hermite unitary polynomial.

\begin{proposition}
	\label{prop:heatequationjacobi}
	Let $m\geq 1$. Assume $p \wedge q > m-1$, then $(t,x) \mapsto \chi_t^{(r,s,m)}(x)$ solves the inverse heat equation:
	\begin{equation*}
		\partial_t\chi_t^{(r,s,m)}(x) = -\left\{\mathcal{L}_x^{(r,s)} + m(r+s+m+1) \right\}
		\chi_t^{(r,s,m)}(x),
	\end{equation*}
	where
	\begin{equation*}
		\mathcal{L}_x^{(r,s)} := x(1-x)\partial_{xx}^2 + [(r+1) - (r+s+2)x]\partial_x
	\end{equation*}
	is the one-dimensional Jacobi operator.
\end{proposition}
\begin{proof}
	Let
	\begin{equation*}
		F(x,\lambda_t^1, \dots, \lambda_t^m) =
		\prod_{i=1}^m(x-\lambda_t^i), \quad x \in \mathbb{R},
	\end{equation*}
	be the characteristic polynomial of the Hermitian Jacobi process $X_t$ and recall that the eigenvalues are almost surely simple. Then, It\^o's formula entails:
	\begin{multline*}
		dF(x, \lambda_t^1, \dots, \lambda_t^m) = \textrm{Local Martingale}
		-F(x, \lambda_t^1, \dots, \lambda_t^m)\sum_{i=1}^m\frac{p-(p+q)\lambda_t^i}{x-\lambda_t^i}
		\\-F(x, \lambda_t^1, \dots, \lambda_t^m)\sum_{i=1}^m\frac{1}{x-\lambda_t^i}
		\sum_{j\neq i}\frac{\lambda_t^i(1-\lambda_t^j) + \lambda_t^j(1-\lambda_t^i)}{\lambda_t^i-\lambda_t^j} dt.
	\end{multline*}
	Note in passing that the quadratic variation terms have zero contribution since
	\begin{equation*}
		\partial_{u_iu_i}F(x, u_1, \dots, u_m) = 0,
	\end{equation*}
	and since the Brownian motions $(B_t^i, t \geq 0)_{1 \leq i \leq m}$ are independent.
	Using
	\begin{equation*}
		\partial_xF(x, u_1, \dots, u_m) = F(x, u_1, \dots, u_m)\sum_{i=1}^m \frac{1}{x-u_i},
	\end{equation*}
	it follows that

	\begin{align*}
		F(x, \lambda_t^1, \dots, \lambda_t^m) \sum_{i=1}^m\frac{p-(p+q)\lambda_t^i}{x-\lambda_t^i}
		 & = p\partial_x F(x,\lambda_t) - (p+q)F(x,\lambda_t)\sum_{i=1}^m \frac{\lambda^i_t}{x-\lambda_t^i} \\
		 & = p\partial_x F(x,\lambda_t) - (p+q)F(x,\lambda_t)\sum_{i=1}^m (\frac{x}{x-\lambda^i_{t}}-1)     \\
		 & = p\partial_x F(x,\lambda_t) - (p+q)(x\partial_x F(x,\lambda_t)-mF(x,\lambda))                   \\
		 & = \partial_x  F(x,\lambda_t)(p - (p+q)x)+m(p+q)F(x,\lambda_t)
	\end{align*}
	where we used the shorthand notation $(x,\lambda_t)$ to denote $(x,\lambda_t^1, \dots, \lambda_t^m)$.
	Furthermore, the term
	\begin{equation*}
		2\sum_{i=1}^m\frac{1}{x-\lambda_t^i} \sum_{j\neq i}\frac{\lambda_t^i(1-\lambda_t^j) + \lambda_t^j(1-\lambda_t^i)}{\lambda_t^i-\lambda_t^j}
	\end{equation*}
	may be symmetrised as
	\begin{align*}
		\sum_{i=1}^m\sum_{j\neq i} \frac{\lambda_t^i(1-\lambda_t^j) + \lambda_t^j(1-\lambda_t^i)}{\lambda_t^i-\lambda_t^j}\left[\frac{1}{x-\lambda_t^i} - \frac{1}{x-\lambda_t^j}\right]
		 & =\sum_{i=1}^m\sum_{j\neq i}\frac{\lambda_t^i(1-\lambda_t^j) + \lambda_t^j(1-\lambda_t^i)}{(x-\lambda_t^i)(x-\lambda_t^j)}
		\\& = 2 \sum_{i=1}^m\sum_{j\neq i}\frac{\lambda_t^i(1-\lambda_t^j)}{(x-\lambda_t^i)(x-\lambda_t^j)}.
	\end{align*}
	On the other hand, one readily computes:
	\begin{equation*}
		\partial_{xx}^2F(x, u_1, \dots, u_m) = F(x, u_1, \dots, u_m)\sum_{i=1}^m\sum_{j \neq i}\frac{1}{(x-u_i)(x-u_j)}
	\end{equation*}
	and split, for any $1 \leq i \neq j \leq m$,
	\begin{align*}
		\frac{x(1-x)}{(x-u_i)(x-u_j)} & = \frac{1-x}{x-u_j} + \frac{u_i(1-x)}{(x-u_i)(x-u_j)}
		\\& = \frac{1-x}{x-u_j} + \frac{u_i(1-u_j)}{(x-u_i)(x-u_j)} - \frac{u_i}{x-u_i}
		\\& = \frac{1-x}{x-u_j} + \frac{u_i(1-u_j)}{(x-u_i)(x-u_j)} - \frac{x}{x-u_i} + 1.
	\end{align*}
	Consequently, again with the shorthand notation $(x,u)$ for $(x,u_1, \dots, u_m)$, one gets:
	\begin{align*}
		x(1-x)\partial_{xx}^2F(x, u) & = [m(m-1) + (m-1)(1-x)\partial_x -
		(m-1)x\partial_x]F(x, u)                                          \\& + F(x, u) \sum_{i=1}^m\sum_{j \neq i} \frac{u_i(1-u_j)}{(x-u_i)(x-u_j)},
	\end{align*}
	whence we deduce that
	\begin{align*}
		F(x, \lambda)  \sum_{i=1}^m & \sum_{j\neq i}\frac{\lambda_t^i(1-\lambda_t^j)}{(x-\lambda_t^i)(x-\lambda_t^j)} = x(1-x)\partial_{xx}^2F(x, \lambda)                                                                                                                        -(m-1)[m+(1-2x)\partial_x]F(x, \lambda).
	\end{align*}
	Gathering all the terms, we get the following SDE for the characteristic polynomial of $X_t$:

	\begin{multline}\label{SDE1}
		d[F(x,\lambda_t)]= \textrm{Local Martingale} +\\
		-x(1-x)\partial_{xx}^2F(x,\lambda_t) - (r+1-(s+r+2)x)\partial_x F(x,\lambda_t)
		- m(d-m+1)F(x,\lambda_t) 
	\end{multline}
	Finally, the bracket of the local martingale part is given by
	\begin{equation*}
		2\sum_{i=1}^m \lambda_t^i(1-\lambda_t^i)\prod_{j \neq i}(x-\lambda_j),
	\end{equation*}
	and is obviously bounded for any fixed $x$ in a bounded interval. Consequently, the local martingale part is a true martingale, and the inverse Jacobi-heat equation follows after taking the expectation of both sides of \eqref{SDE1}.
\end{proof}

It is known that the Jacobi operator $\mathcal{L}_x^{(r,s)}$ admits a complete set of eigenpolynomials given by the orthogonal Jacobi polynomials $Q_j^{(r,s)}(x), j \geq 0$ in $[0,1]$:
\begin{equation*}
	Q_j^{(r,s)}(x) := \frac{(r+1)_j}{j!}{}_2F_1(-j, r+s+j+1, r+1; x),
\end{equation*}
where ${}_2F_1$ stands for the Gauss hypergeometric function. More precisely, one has:
\begin{equation*}
	\mathcal{L}_x^{(r,s)}[Q_j^{(r,s)}](x)= -j(j+r+s+1)Q_j^{(r,s)}(x).
\end{equation*}
Consequently,
\begin{align*}
	\left\{\mathcal{L}_x^{(r,s)} + m(d-m+1)\right\}Q_j^{(r,s)}(x) & =
	\left\{-j(j+r+s+1) + m(r+s+m+1)\right\}Q_j^{(r,s)}(x)
	\\& = (m-j)(r+s+1+m+j)Q_j^{(r,s)}(x).
\end{align*}
\subsection{Proof of Theorem \ref{cor:expressionaveragedcharacteristic}}
It is known that the Jacobi operator $\mathcal{L}_x^{(r,s)}$ admits a complete set of eigenpolynomials, the orthogonal Jacobi polynomials :
\begin{equation*}
	\mathcal{L}_x^{(r,s)}[Q_j^{(r,s)}](x)= -j(j+r+s+1)Q_j^{(r,s)}(x).
\end{equation*}
Consequently,
\begin{align*}
	\left\{\mathcal{L}_x^{(r,s)} + m(d-m+1)\right\}Q_j^{(r,s)}(x) & =
	\left\{-j(j+r+s+1) + m(r+s+m+1)\right\}Q_j^{(r,s)}(x)
	\\& = (m-j)(r+s+1+m+j)Q_j^{(r,s)}(x).
\end{align*}

Proposition \ref{prop:heatequationjacobi} then yields the first part of Theorem \ref{thm:maintheorem}.
It suffices to find the expansion of $(x-1)^m$ in the Jacobi polynomial basis. But this is afforded by the following formula \cite{MR1697414}:
	\begin{equation*}
		\left(\frac{1-u}{2}\right)^m = \sum_{j=0}^m \frac{(-m)_j(\alpha+j+1)_{m-j}(\alpha+\beta+2j+1)}{(\alpha+\beta+j+1)_{m+1}}P_j^{(\alpha,\beta)}(u), \quad u \in [-1,1],
	\end{equation*}
	where $P_j^{(\alpha,\beta)}$ is the $j$-th Jacobi polynomial in $[-1,1]$:
	\begin{equation*}
		P_j^{(\alpha,\beta)}(1-2x) = Q_j^{(\alpha,\beta)}(x), \quad \alpha, \beta > -1.
	\end{equation*}
	Indeed, the symmetry property $P_j^{(\alpha,\beta)}(u) = P_j^{(\beta, \alpha)}(-u)$ transforms the last formula into
	\begin{equation*}
		\left(\frac{1+u}{2}\right)^m = \sum_{j=0}^m (-1)^j\frac{(-m)_j(\alpha+j+1)_{m-j}(\alpha+\beta+2j+1)}{(\alpha+\beta+j+1)_{m+1}}P_j^{(\beta,\alpha)}(u).
	\end{equation*}
	Substituting $u = 1-2x, \beta = r, \alpha = s$ there, we obtain the expansion:
	\begin{align}
    \label{eqn:finiteexpansion}
		(1-x)^m = m!\sum_{j=0}^m \frac{\Gamma(s+m+1)\Gamma(r+s+j+1)(r+s+2j+1)}
		{(m-j)!\Gamma(j+s+1)\Gamma(r+s+m+j+2)}Q_j^{(r,s)}(x)
	\end{align}
	where we used the representation of the Pochhammer symbol valid for positive real numbers: 
    \begin{equation*}
    (y)_j = \frac{\Gamma(y+j)}{\Gamma(y)}, \quad y > 0. 
    \end{equation*}
Applying the shifted Jacobi-heat semi-group
$   e^{-t(\mathcal{L}_x^{(r,s)}+m(d-m+1))}
$
to \eqref{eqn:finiteexpansion}, we are done.
\begin{remark}
\label{rk:longtime}
	As $t \rightarrow +\infty$, only the term $j=m$ in the above expansion of 
    $\chi_t^{(r,s,m)}(x)$ gives a non zero contribution.
	Accordingly, the limiting averaged characteristic polynomial
	\begin{equation*}
		\chi_{\infty}^{(r,s,m)}(x) := \lim_{t \rightarrow +\infty} \chi_t^{(r,s,m)}(x)
	\end{equation*}
	is proportional to the monic Jacobi polynomial $Q_{m}^{(r,s)}(x)/k_m^{(r,s)}$, where $k_m^{(r,s)}$ is the leading term of $Q_{m}^{(r,s)}(x)$. This is in agreement with the fact that the averaged characteristic polynomial of the Jacobi unitary ensemble, the weak limit of the Hermitian Jacobi process as $t \rightarrow +\infty$, is given by a Jacobi polynomial. Indeed, the eigenvalues of the JUE (Jacobi Unitary Ensemble) are given by a Selberg weight, as such its averaged characteristic polynomial is an instance of the celebrated Aomoto integral (see e.g. \cite{MR1985318}).
\end{remark}
For any fixed time $t \geq 0$, recall
\begin{equation}
	\label{eqn:averagedpolynomial}
	\chi_t^{(r,s,m)}(x) = \prod_{m \geq k \geq 1}(x-x_{k}^{(r,s,m)}(t))= \mathbb{E}\left[\prod_{m \geq i \geq 1}(x-\lambda_i(t))\right].
\end{equation}
\begin{proposition}
	\label{prop:positivityoftheroot}
	The roots $(x_j^{(r,s,m)}(t))_{1 \leq j\leq m}$ of $\chi_{t}^{(r,s,m)}$ are all real. In addition, up to re-indexing, they satisfy the following ODE:
	\begin{align}\label{ODEroot}
		\frac{dx_j^{(r,s,m)}}{dt}(t) & = (p - (p+q)x_j^{(r,s,m)}(t)) \\&\hspace{2cm}+
		\sum_{k \neq j} \frac{x_j^{(r,s,m)}(t)(1-x_k^{(r,s,m)}(t))+x_k^{(r,s,m)}(t)(1-x_j^{(r,s,m)}(t))}{x_j^{(r,s,m)}(t)-x_k^{(r,s,m)}(t)},		\nonumber
		\\&  = (r+1) - (r+s+2)x_j^{(r,s,m)}(t) 		\nonumber
        \\ &\hspace{2cm}+ 2x_j^{(r,s,m)}(t)(1-x_j^{(r,s,m)}(t))
		\sum_{k \neq j} \frac{1}{x_j^{(r,s,m)}(t)-x_k^{(r,s,m)}(t)}.
		\nonumber
	\end{align}
	Besides, for any time $t > 0$ and any $1 \leq j \leq m$, $x_j^{(r,s,m)}(t) \in (0,1)$.
\end{proposition}
\begin{proof} The first part of this proposition is due to Voit, as stated in Theorem 4.4 of \cite{MR4781073}. Now, let us turn to the second statement. We have
\begin{equation*}
	\prod_{k=1}^mx_k^{(r,s,m)}(t) = \mathbb{E}\left[\prod_{k=1}^m\lambda_k(t)\right] > 0
\end{equation*}
since $\lambda_k(t) > 0$ almost surely for any $1 \leq k \leq m$ \cite{MR2719370, MR3296535}. It follows that $x_k^{(r,s,m)}(t) \neq 0$ for any $1 \leq k\leq m$. Similarly,
\begin{equation*}
	e_k^{(r,s,m)}(x_1^{(r,s,m)}(t), \dots, x_m^{(r,s,m)}(t)) = \mathbb{E}[e_k^{(r,s,m)}(\lambda_1(t), \dots, \lambda_m(t))] > 0
\end{equation*}
which implies that $x_k^{(r,s,m)}(t) > 0$ for any $1 \leq k\leq m$. In fact, if all the elementary symmetric polynomials on real numbers $x_1, \dots, x_m$ are positive, 
then so do $(x_1, \dots, x_m)$. This is proved by setting
\begin{equation*}
	J(x):= \prod_{k=1}^m(x+x_k)
\end{equation*}
and argue that if $x_i < 0$ for some $i \in [m]$, this yields the following contradiction:
\begin{equation*}
	0 = J(-x_i) = \sum_{k=1}^m (-x_i)^{m-k}e_k(x_1,\dots, x_m) > 0.
\end{equation*}
To prove now that $x_k^{(r,s,m)}(t) < 1$ for any $1 \leq k \leq m$, observe that 
\begin{equation*}
	\prod_{k=1}^m (x-(1-x_k^{(r,s,m)}(t)) =  \prod_{k=1}^m (x_k^{(r,s,m)}(t)-(1-x))
	= \mathbb{E}\left[\prod_{k=1}^m (x-(1-\lambda_k(t))\right].
\end{equation*}
But $(1-\lambda_k(t), t \geq 0)_{1 \leq k \leq m}$ is a Jacobi particle system  with parameters $(s,r,m)$. As a matter of fact, $ 1-x_k^{(r,s,m)}(t) > 0$.
\end{proof}

\subsection{Relation to the Hermite unitary polynomial}
\label{ssec:unitaryhermite}
In \cite{MR3071702}, it was shown that the spectral distribution of the free Jacobi process with parameters $(1,1/2)$ in the compressed probability space coincides with the spectral distribution of 
\begin{equation*}
 \frac{1}{4}(Y_{2t} + Y^{\star}_{2t} + 2) = \frac{1}{2}(1+\Re(Y_{2t})),    
\end{equation*}
where $(Y_t)_{t \geq 0}$ is the free unitary Brownian motion \cite{MR1430721}. Note in this respect that this special set of parameters corresponds at the matrix level to Hermitian Jacobi processes, which are radial parts of `asymptotically square' (as $  m\rightarrow +\infty$) principal minors of a unitary Brownian motion whose size is asymptotically twice the sizes of these truncations.    

On the other hand, the finite analogue of the free unitary convolution, with respect to which $(Y_t)_{t \geq 0}$ is a free L\'evy process, was introduced and studied in \cite{MR4285332}. 
This finite convolution is encoded by the zeroes of the Hermite unitary polynomial defined for any time $t \geq 0$ by \cite{MR4912666}\footnote{This definition differs from the original one given in \cite{MR4285332} by the time change $t \mapsto -t/(d-1)$.}: 
$$
H_d(z,t) = \sum_{k=0}^d z^{d-k} (-1)^k \binom{d}{k} \exp\left(-t\frac{k(d-k)}{2}\right).
$$
In particular, Lemma 2.1. in \cite{MR4912666} asserts that the zeroes of $H_d(\cdot, t)$ lie on the unit circle and we infer from Corollary 3.35 proved in \cite{MR4285332} that the dynamics of the corresponding angles satisfy the following ODE: 
\begin{equation}\label{AngularDyn}
\partial_t \theta_j(t) = \frac{1}{2}\sum_{k\neq j} \cot\left(\frac{\theta_j(t) - \theta_k(t)}{2}\right), \quad 1 \leq j \leq d.
\end{equation}
In particular, Theorem 6.1 in \cite{MR4781073} shows that these angles (and in turn the roots of $H_d(\cdot, t)$) remain distinct for all times $t > 0$ even if they collapse at $t=0$. 
Besides, since $H_d(\cdot, t)$ has real coefficients, then its roots are pairwise-conjugate. Note also that $z=-1$ is not a root of $H_d(\cdot, t)$ since $H_d(-1, t)$ is a sum of positive terms, 
while $z=1$ is so only when the degree $d$ is odd. 
Based on this discussion, it is tempting to wonder whether both root dynamics \eqref{AngularDyn} and \eqref{ODEroot} are related through the transformation\footnote{For sake of simplicity, we omit the dependence on the parameters.}: 
\begin{equation*}
x_j(t) = \frac{1+\cos(\theta_j(2t))}{2} = \cos^{2}\left(\frac{\theta_j(2t)}{2}\right), \quad 1 \leq j \leq m,  
\end{equation*}
when $d = 2m$ is even, where here $0< \theta_1(t) < \theta_2 < \cdots < \theta_m < \pi$ are the ordered angles of the zeroes of $H_d$ lying in the upper-half of the unit circle. 
The following proposition shows that this is indeed the case: 
\begin{proposition} Assume $\chi^{(-1/2,-1/2,m)}_0(x) = (x-1)^{2m} = H_{2m}(x,0)$. Then, for any $t\geq 0$ and any $m\geq 0$:
\label{prop:relationunitaryjacobi}
    \begin{align}
    H_{2m}(z,2t)=4^mz^m\chi^{(-1/2,-1/2,m)}_t\left(\frac{z+z^{-1}+2}{4}\right).
    \end{align}
\end{proposition}
    
\begin{proof}
For $m=0$, the result is immediate. 
We turn our attention to the case $ m\ge 1$. Set $\eta_j(t) = \cos^2(\theta_j(t)/2)$, for $1 \leq j \leq m$. Then: 
\begin{align*}
\partial_t \eta_j(t) &= -\left[\partial_t \theta_j(t)\right] \sin(\theta_j(t)/2)\cos(\theta_j(t)/2)\\
&= -\frac{1}{2}\sum_{\substack{k\neq j \\ 1 \leq k \leq 2m}} \cot\left(\frac{\theta_j(t) - \theta_k(t)}{2}\right) \sin(\theta_j(t)/2)\cos(\theta_j(t)/2) \\
&= -\frac{1}{2}\sum_{\substack{k\neq j \\ 1 \leq k \leq 2m}} \frac{\cos^2(\theta_j(t)/2)\cos(\theta_k(t)/2)\sin(\theta_j(t)/2) 
+ \sin^2(\theta_j(t)/2)\sin(\theta_k(t)/2)\cos(\theta_j(t)/2) }{\sin((\theta_j(t) - \theta_k(t))/2)}.
\end{align*}
Using the identity $\sin^2(\theta_j(t)/2) = 1- \cos^2(\theta_j(t)/2)$, we further get:
\begin{align*}
\partial_t \eta_j(t) &= -\frac{2m-1}{2}\eta_j(t) - \frac{1}{2} \sum_{\substack{k\neq j \\ 1 \leq k \leq 2m}} \frac{\sin(\theta_k(t)/2)\cos(\theta_j(t)/2)}{\sin((\theta_j(t) - \theta_k(t))/2)} \\
&= - \frac{2m-1}{2}\eta_j(t) - \frac{1}{2} \sum_{\substack{k\neq j \\ 1 \leq k \leq 2m}}\left(\frac{\tan(\theta_j(t)/2)}{\tan(\theta_k(t)/2)} - 1\right)^{-1}.
\end{align*}
Taking out the summand corresponding to $k=2m-j+1$ and remembering that the angles come into opposite pairs, the last sum splits into:
\begin{align*}
\sum_{\substack{k\neq j \\ 1 \leq k \leq 2m}}\left(\frac{\tan(\theta_j(t)/2)}{\tan(\theta_k(t)/2)} - 1\right)^{-1} & = -\frac{1}{2} +
\sum_{\substack{k\neq j \\ 1 \leq k \leq m}}\left(\frac{\tan(\theta_j(t)/2)}{\tan(\theta_k(t)/2)} - 1\right)^{-1}  
 - \left(\frac{\tan(\theta_j(t)/2)}{\tan(\theta_k(t)/2)} + 1\right)^{-1} \\
 \\
&= -\frac{1}{2} + 2\sum_{\substack{k\neq j \\ 1 \leq k \leq m}}\frac{\tan^2(\theta_k(t)/2)}{{\tan^2(\theta_j(t)/2)}-{\tan^2(\theta_k(t)/2)}} \\
&= -\frac{1}{2}+ 2\sum_{\substack{k\neq j \\ 1 \leq k \leq m}}\frac{\cos^2(\theta_j(t)/2)-\cos^2(\theta_j(t)/2)\cos^2(\theta_k(t)/2)}{{\cos^2(\theta_k(t)/2)}-{\cos^2(\theta_j(t)/2)}}.
\end{align*}
As result, for any $1 \leq j \leq m$, the map $\eta_j$ satisfies the ODE:
\begin{align*}
\partial_t \eta_j(t) &= -\frac{2m-1}{2}\eta_j(t) + \frac{1}{4} -  \sum_{\substack{k\neq j \\ 1 \leq k \leq m}}
\frac{\cos^2(\theta_j(t)/2)-\cos^2(\theta_j(t)/2)\cos^2(\theta_k(t)/2)}{\cos^2(\theta_k(t)/2)-\cos^2(\theta_j(t)/2)} \\
&= -\frac{2m-1}{2}\eta_j(t) + \frac{1}{4}  + \sum_{\substack{k\neq j \\ 1 \leq k \leq m}} \frac{\eta_j(t)(1-\eta_k(t))}{\eta_j(t)-\eta_k(t)}
\end{align*}
Equivalently, $\tilde{\eta}_j(t) := \eta_j(2t)$ satisfies: 
\begin{align*}
\partial_t \tilde{\eta_j}(t) 
&=-(2m-1)\tilde{\eta}_j(t)  + \frac{1}{2} + 2 \sum_{\substack{k\neq j \\ 1 \leq k \leq m}} \frac{\tilde{\eta_j}(t)(1-\tilde{\eta}_k(t))}{\tilde{\eta}_j(t)-\tilde{\eta}_k(t)}
\end{align*}
Finally, notice that 
\begin{align*}
 2\frac{\tilde{\eta}_j(t)(1-\tilde{\eta}_k(t))}{\tilde{\eta}_j(t)-\tilde{\eta}_k(t)} &= 
 \frac{\tilde{\eta}_j(t)(1-\tilde{\eta}_k(t)) + \tilde{\eta}_k(t)(1-\tilde{\eta}_j(t))}{\tilde{\eta}_j(t)-\tilde{\eta}_k(t)} + 
 \frac{\tilde{\eta}_j(t)(1-\tilde{\eta}_k(t)) - \tilde{\eta}_k(t)(1-\tilde{\eta}_j(t))}{\tilde{\eta}_j(t)-\tilde{\eta}_k(t)} \\
 &= 1+ \frac{\tilde{\eta}_j(t)(1-\tilde{\eta}_k(t)) + \tilde{\eta}_k(t)(1-\tilde{\eta}_j(t))}{\tilde{\eta}_j(t)-\tilde{\eta}_k(t)},
\end{align*}
we end up with:
\begin{align*}
\partial_t \tilde{\eta_j}(t) 
&= \left(m-\frac{1}{2}\right)  -(2m-1)\tilde{\eta}_j(t) +
 \sum_{\substack{k\neq j \\ 1 \leq k \leq m}}\frac{\tilde{\eta}_j(t)(1-\tilde{\eta}_k(t)) + \tilde{\eta}_k(t)(1-\tilde{\eta}_j(t))}{\tilde{\eta}_j(t)-\tilde{\eta}_k(t)}. 
\end{align*}
This is the ODE \eqref{ODEroot} with $p=q = m- (1/2)$. But Theorem 1.1. in \cite{MR4751304} shows (after performing an affine variable change) 
that \eqref{ODEroot} admits a unique solution for any $p,q > m-1$ which remains in the open domain $\{0 < x_m < \cdots < x_1 < 1\}$ at any time $t > 0$. It follows that 
\begin{equation*}
\cos^2\left(\frac{\theta_j(2t)}{2}\right) = \tilde{\eta_j}(t) = x_j(t)    
\end{equation*}
for any $1 \leq j \leq m$ and any $t \geq 0$, therefore 
\begin{align*}
H_{2m}(z,2t) & = \prod_{j=1}^m|z-e^{i\theta_j(2t)}|^2 
\\& = \prod_{j=1}^m(z^2-2z \cos(\theta_j(2t)) + 1)
\\& = \prod_{j=1}^m(z^2-2z (2x_j(t) - 1)+1) = (4z)^m \chi^{(-1/2,-1/2,m)}_t\left(\frac{z+z^{-1}+2}{4}\right).
\end{align*}
\end{proof} 
\begin{remark}
    We can infer a relation between the Hermitian and Jacobi heat generators. Let 
    $\mathbb{R}_{[0,1]}[x]_m$ be the set of monic polynomials with degree $m$ and roots in $[0,1]$ and define:
    $$
\varphi_{m}\colon \mathbb{R}_{[0,1]}[x]_m \to \mathbb{R}[z]_{2m}\quad P \mapsto 4^mz^mP(\frac{1}{4}(z+z^{-1}+2)).
    $$
    Recall the definition of the Hermitian heat generators $\mathcal{L}^{(2m)}_{A}$ and the relation with the Hermitian Unitary polynomial:
    $$
H_{2m}(z,2t)=\exp(-t(z\partial_z){(2m-z\partial_z)})\{(z-1)^{2m}\} = \exp(\mathcal{L}_{A}^{2m})\{(z-1)^{2m} \}.
$$
Then,
    \begin{equation}
    \label{eqn:generators}
    \varphi_m\circ \mathcal{L}^{(-1/2,-1/2,m)} = \mathcal{L}_{A}^{(2m)} \circ \varphi_m.
    \end{equation}

The polynomial $
\varphi_m(P)
$
has real coefficients, hence its roots are pairwise complex-conjugate. Moreover, $z^{2m}\varphi_m(\frac{1}{z})=\varphi_m(z)$ (it is invariant by inversion). Hence, the set of roots of $
\varphi_m(P)
$ is invariant under inversion and conjugation. Moreover, if $z_0$ of $\varphi_m(P)$
then $ \Im(z_0 + z_0^{-1}) = 0$ since $\frac{1}{4}(z_0 + z_0^{-1}+2)$ is root of $P$. Hence, $|z_0| -|z_0|^{-1} = 0$ and $|z_0| = 1$. Thus $\varphi_m(P)$ has roots on the unit circle and they are pair-wise conjugate.
The dynamic of the angles of the roots of $\exp(t\mathcal{L}_{A}^{(2m)})\{\varphi_m(P)\}$ is prescribed by \eqref{AngularDyn}. Equality \eqref{eqn:generators} follows from the same reasoning as exposed in the proof of Proposition \ref{prop:relationunitaryjacobi}.

\end{remark}

\subsection{Another formula for $\chi_t^{(r,s,m)}$}
\label{ssec:anotherformula}
In this paragraph, We derive another formula for $\chi_t^{(r,s,m)}$ based on the explicit expression of the heat kernel of the eigenvalue process $(\lambda_t^{(r,s,m)})_{t \geq 0}$ as well as on 
an identity arising in the study of the infinite unitary group. 
Actually, this kernel admits a series representation in terms of multivariate (symmetric) Jacobi polynomials which are mutually orthogonal with respect to the unitary Selberg weight. 
In this regard, recall that this property is not satisfied by any orthogonal set of multivariate polynomials since it requires the orthogonality of any two polynomials corresponding to different partitions. Moreover, an explicit expansion of the characteristic polynomial of $X_t(m)$ in this polynomial basis may be deduced from the so-called dual Cauchy formula proved in \cite{MR3498194}.  
In particular, the derivation of the obtained formula does not use the (deterministic) dynamics of the roots of the $\chi_t^{(r,s,m)}$.

 For the sake of completeness, we provide a brief reminder of the key facts and results necessary to prove the expansion below. We refer the reader to the paper \cite{MR2719370} for more details.

Let 
\begin{equation*}
	\tau=(\tau_1 \geq \tau_2 \geq ... \geq \tau_m \geq 0)
\end{equation*}
be a partition of length at most $m$ and let $(\tilde{Q}_j^{(r,s)})_{j \geq 0}$ be the sequence of orthonormal Jacobi polynomials with respect to the beta weight:
\begin{equation*}
	u^r(1-u)^s {\bf 1}_{[0,1]}(u).
\end{equation*}
These are given by
\begin{equation*}
	\tilde{Q}_j^{(r,s)}(x) := \frac{Q_j^{(r,s)}}{||Q_j^{(r,s)}||_2} = \left[\frac{(2j+r+s+1)\Gamma(j+r+s+1)j!}{\Gamma(r+j+1)\Gamma(s+j+1)}\right]^{1/2}Q_j^{(r,s)}(x).
\end{equation*}
Then the orthonormal multivariate Jacobi polynomial corresponding to $\tau$ is defined by:
\begin{align*}
	\tilde{Q}_{\tau}^{(r,s,m)}(y_1,...,y_m) := \frac{\det(\tilde{Q}_{\tau_i-i+m}^{(r,s)}(y_j))_{1\leq i,j \leq m}}{V(y_1,...,y_m)},
\end{align*}
where
\begin{equation*}
	V(y_1,...,y_m) := \prod_{1\leq i<j\leq m}(y_i-y_j),
\end{equation*}
is the Vandermonde determinant. If the coordinates $(w_1,...,w_m)$ overlap, this definition still makes sense by either applying L'H\^opital's rule or equivalently by using the expansion of $\tilde{Q}_{\tau}^{(r,s,m)}$ in the Schur polynomial basis.

The multivariable Jacobi polynomials are symmetric (invariant under permutations) and satisfy the remarkable property of being mutually orthonormal with respect to the unitary Selberg weight:
\begin{equation*}
	W^{(r,s,m)}(y_1,\dots, y_m):= [V(y_1, \dots, y_m)]^2\prod_{i=1}^m y_i^r(1-y_i)^s{\bf 1}_{0<y_m<...<y_1<1}, \quad r,s > -1.
\end{equation*}
More precisely, for any two different partitions $\tau$ and $\kappa$, one has:
\begin{equation*}
	\int Q_{\tau}^{(r,s,m)}(y_1,...,y_m)Q_{\kappa}^{(r,s,m)}(y_1,...,y_m)W^{(r,s,m)}(y_1,\dots, y_m)dy_1\cdots dy_m = 0,
\end{equation*}
as one readily checks using Andreief's identity. 

Replacing $\tilde{Q}_j^{(r,s)}$ with $Q_j^{(r,s)}$, we get the orthogonal multivariate Jacobi polynomials $(Q_{\tau}^{(r,s,m)})_{\tau}$ in $[0,1]^m$, and performing further the variable change $y \mapsto 1-2y$, one gets the orthogonal multivariate Jacobi polynomials $(P_{\tau}^{(r,s,m)})_{\tau}$ in $[-1,1]^m$ and used in \cite{MR3498194}.

Now, recall that the semi-group density of the eigenvalues process of $X_t(m)$ and starting at $w$ admits the following absolutely-convergent expansion:
\begin{align*}
	G_t^{r,s,m}(w,y) := \sum_{\tau} e^{-\nu_\tau t}\tilde{Q}_\tau^{(r,s,m)}(w)\tilde{Q}_\tau^{r,s,m}(y) W^{r,s,m}(y_1,\dots, y_m), \quad r,s > -1,
\end{align*}
where
\begin{align*}
	\nu_{\tau} :=\sum_{i=1}^m\tau_i(\tau_i+r+s+1+2(m-i)).
\end{align*}

Now, we are ready to prove:
\begin{proposition}\label{Generalexpansion}
	For any $w \in [0,1]^m$, the following expansion holds:
	\begin{equation*}
		\chi_t^{(r,s,m)}(x) = \frac{1}{(-2)^{m(m+1)/2}}  \sum_{j=0}^m (-1)^{m-j}  e^{-\nu_{(1^{m-j})} t}  \frac{Q_{1^{m-j}}^{(r,s,m)}(w)}{k_{1^{m-j}}^{(r,s,m)}} \frac{Q_{j}^{(r,s)}(x)}{k_{j}^{(r,s)}},
	\end{equation*}
	where \begin{itemize}
	    \item $1^{m-j}$ is the partition with only $(m-j)$ ones (if $j=m$ then it is the zero partition).
        \item $k_{1^{m-j}}^{(r,s,m)}$ is the leading coefficient of $P_{1^{m-j}}^{(r,s,m)}$ when expanded in the basis of Schur polynomials (\cite{MR3075094}, eq. (3)). 
        \item $k_{j}^{(r,s)} := k_{j}^{(r,s,1)}$ is the leading coefficient of $P_j^{(r,s)} := P_{j}^{(r,s,1)}$.
    \end{itemize}
\end{proposition}
\begin{proof}
    Recall from \cite{MR3075094} (p. 272) the dual Cauchy-identity for the multivariate Jacobi polynomials in $[-1,1]^m$:
	\begin{equation*}
		\prod_{i=1}^N\prod_{j=1}^K (u_i + v_j) = \sum_{\substack{\lambda = (\lambda_1 \geq \dots \lambda_N \geq 0) \\ \lambda_1 \leq K}} \frac{P_{\mu}^{(s,r,m)}(v_1,\dots, v_K) P_{\lambda}^{(r,s,m)}(u_1, \dots, u_N)}{k_{\mu}^{(s,r,m)}k_{\lambda}^{(r,s,m)}},
		\quad u_i, v_j \in [-1,1],
	\end{equation*}
	where for any partition $K \geq \lambda_1 \geq \lambda_2 \geq \cdots \geq \lambda_N$, the partition $\mu$ is defined through the conjugate partition
    $\lambda' = (\lambda'_1 \geq \dots \geq \lambda'_K)$ of $\lambda$ by:
	\begin{equation*}
		\mu = (N - \lambda'_K, \dots, N - \lambda'_1).
	\end{equation*}
 
 Now, the hypergeometric representation 
    \begin{equation*}
        P_j^{(r,s)}(u) = \frac{(r+1)_j}{j!}{}_2F_1\left(-j, r+s+j+1, r+1; \frac{1-u}{2}\right)
    \end{equation*}
shows that 
\begin{equation*}
  k_{j}^{(r,s,1)} = \frac{(r+s+1+j)_j}{j!2^j}.  
\end{equation*}
In turn, the determinantal form of $P_j^{(r,s,m)}$ entails 
\begin{equation*}
  k_{\lambda}^{(r,s,m)} = \prod_{j=1}^m k_{\lambda_j+m-j}^{(r,s,1)},  
\end{equation*}
for any partition $\lambda$. 

	We specialise the dual Cauchy formula by putting 
    $$
    N=1,\quad K= m, \quad u_1 = 1-2x, \quad v_j = 2y_j-1.
    $$ 
    Then $\lambda = \lambda_1 \in \{0, \dots, m\}$ is a row so that $\lambda' = 1^{\lambda_1}$ is a column and in turn $\mu = 1^{m-\lambda_1}$ is a column as well.
	Using the mirror property satisfied by the Jacobi polynomials in $[-1,1]$, we get
	\begin{align*}
		(-2)^mF(x, y_1, \dots, y_m) & = \sum_{j=0}^m \frac{P_{1^{m-j}}^{(s,r,m)}(2y_1-1,\dots,
		2y_m-1) P_{j}^{(r,s,1)}(1-2x)}{k_{1^{m-j}}^{(r,s,m)}k_{j}^{(r,s,1)}}
		\\& = \sum_{j=0}^m (-1)^{m-j} \frac{P_{1^{m-j}}^{(r,s,m)}(1-2y_1,\dots,
		1-2y_m) P_{j}^{(r,s)}(1-2x)}{k_{1^{m-j}}^{(r,s,m)}k_{j}^{(r,s,1)}}
		\\& = \sum_{j=0}^m (-1)^{m-j} \frac{Q_{1^{m-j}}^{(r,s,m)}(y_1,\dots,
		y_m) Q_{j}^{(r,s)}(x)}{(-2)^{m(m-1)/2}k_{1^{m-j}}^{(r,s,m)}k_{j}^{(r,s,1)}}.
	\end{align*}
	Now, we can write the heat kernel $G_t^{(r,s,m)}$ as:
	\begin{align*}
		G_t^{(r,s,m)}(w,y) := \sum_{\tau} e^{-\nu_\tau t}
		\frac{Q_\tau^{(r,s,m)}(w)Q_\tau^{(r,s,m)}(y)}{\left(\prod_{j=1}^m||Q_{\tau_j+m-j}^{(r,s)}||_2\right)^2} W^{(r,s,m)}(y_1,\dots, y_m),
	\end{align*}
	then appeal to the mutual orthogonality of $(Q_\tau^{(r,s,m)})_{\tau}$ to compute the integral:
	\begin{equation*}
		\chi_t^{(r,s,m)}(x) = \int \prod_{i=1}^m(x-y_i)G_t^{(r,s,m)}(w,y) dy_1\cdots dy_m.
	\end{equation*}
	Doing so only leaves the partitions $\tau = 1^{m-j}, 0 \leq j \leq m$, whence the expansion follows. 
    \end{proof}
\begin{remark}
	Proposition 7.1 in \cite{MR3075094} gives an explicit expression of
	$Q_\tau^{(r,s,m)}(1^m)$ as a ratio of products of Gamma functions. After lengthy (but easy) computations, one retrieves the second statement of Corollary \ref{cor:expressionaveragedcharacteristic} and corresponds to the starting point $w = 1^m$. 
    
    On the other hand, if $(w_i = z_i^{(r,s)})_{1 \leq i \leq m}$ are the zeroes of the Jacobi polynomials $Q_m^{(r,s)}$, then $Q_{(1^{m-j})}^{(r,s,m)}(w) = 0$ for all $0 \leq j \leq m-1$ since $Q_{\tau_{1+m-1}}^{(r,s)}(z_i) = 0$ for any $1 \leq i \leq m$. Consequently, $\chi_t^{(r,s,m)}(x) = Q_m^{(r,s)}(x)/k_m^{(r,s)}$ for any $t \geq 0$ and agrees with the fact that $(z_i)_{1 \leq i \leq m}$ is the stationary solution of \eqref{ODEroot} (see Proposition 7 in \cite{MR4673381}).
\end{remark}
\section{Frozen Hermitian Jacobi process and finite free probability} \label{sec:frozen}

In this section, we study the Frozen Jacobi process $(x^{(r,s,m)}_t)_{t \geq 0}$ in the limit as $m$ tends to infinity and its finite free $T$-transform. We prove Theorem \ref{thm:maintheorem} as a corollary of Proposition \ref{prop:momentsroots}, where we derive the ordinary differential equation satisfied by the moment sequence of this counting measure, using a relative compactness argument. Afterwards, we proceed to the proof of Theorem \ref{thm:finitefreedifference} where we distinguish three cases according to the location of the roots of the considered polynomials. 
Next, Corollary \ref{cor:odefiniteTtransform} encloses a differential-difference equation for the finite free $T$-transform. Passing to the limit there, we retrieve in Theorem \ref{thm:PDEfromfinitefree} the second PDE displayed in Theorem \ref{thm:derivationfreestransform} and satisfied by the free $T$-transform. 

\subsection{Proof of Theorem \ref{thm:maintheorem}}
For any time $t \geq 0$, let
\begin{equation*}
	m_{\ell}^{(r,s,m)}(t) =  \mu^{(r,s,m)}_t(x^\ell), \quad \ell \geq 0,
\end{equation*}
be the moment sequence of $\mu_t^{(r,s,m)}$.
\begin{proposition}\label{prop:momentsroots}
	For any $r,s \geq 0$ and any integer $m \geq 1$, the trajectory of the infinite sequence of moments $(m_{\ell}^{(r,s,m)}(t))_{\ell\geq 0}$ is the solution of the following lower triangular differential system :
	\begin{equation*}
		\frac{d}{dt}m_1^{(r,s,m)}(t)=p  -(p+q) m_1^{(r,s,m)}(t),
	\end{equation*}
	and for any $\ell\ge2$ :
	\begin{align}\label{ODEmoments}
		\frac{d}{dt}m_{\ell}^{(r,s,m)}(t)&= -\ell(p+q-\ell+1) m_{\ell}^{(r,s,m)}(t) \\
        &\nonumber \hspace{1cm}+ \ell(p-\ell+1)  m_{\ell-1}^{(r,s,m)}(t)
		+m\ell\sum_{k=0}^{\ell-2}(m_{k}^{(r,s,m)}(t)-m_{k+1}^{(r,s,m)}(t)m_{\ell-1-k}^{(r,s,m)}(t).
	\end{align}
\end{proposition}
\begin{proof}
	Appealing to the ODE \eqref{ODEroot}, we derive
    \begin{align*}
		\frac{d}{dt}m_{\ell}^{(r,s,m)}(t) & =\frac{\ell}{m}\sum_{i=1}^m (x^i(t))^{\ell-1}\frac{d}{dt}x^i(t)
		\\&=\frac{\ell}{m}\bigg(mp  m_{\ell-1}^{(r,s,m)}(t)-m(p+q) m_{\ell}^{(r,s,m)}(t)(t)
        \\& \hspace{2cm}+\sum_{i=1}^m\sum_{j\neq i} \frac{(x^i(t))^\ell(1-x^j(t)) + (x^i(t))^{\ell-1}x^j(t)(1 - x^i(t))}{x^i(t) - x^j(t)}\bigg).
	\end{align*}
	Now, we expand
	\begin{align*}
		&\sum_{i=1}^m\sum_{j\neq i} \frac{(x^i(t))^\ell(1-x^j(t)) + (x^i(t))^{\ell-1}x^j(t)(1 - x^i(t))}{x^i(t) - x^j(t)} \\
        &\hspace{5cm}=\sum_{i=1}^m\sum_{j\neq i} \frac{-2(x^i(t))^\ell x^j(t) }{x^i(t) - x^j(t)}+\frac{(x^i(t))^\ell+(x^i(t))^{\ell-1}x^j(t)}{x^i(t)- x^j(t)},
	\end{align*}
	and make the first double sum for $l\ge 2$ symmetric as:
	\begin{align*}
		-2\sum_{i=1}^m\sum_{j\neq i} \frac{(x^i(t))^{\ell-1}(x^i(t)x^j(t)) }{x^i_t - x^j(t)} & =-2\sum_{i=1}^m\sum_{j\neq i} \frac{x^i(t)x^j(t)[(x^i(t))^{\ell-1}-(x^j(t))^{\ell-1}] }{x^i_t - x^j(t)}
		\\&=-\sum_{k=0}^{\ell-2}\sum_{i=1}^m\sum_{j\neq i} x^i(t)x^j(t)(x^i(t))^k(x^j(t))^{\ell-2-k}
		\\&=-\sum_{k=0}^{\ell-2}\sum_{i=1}^m\sum_{j\neq i} (x^i(t))^{k+1}(x^j(t))^{\ell-1-k}
		\\&=-m^2\sum_{k=0}^{\ell-2}m_{k+1}^{(r,s,m)}(t)m_{\ell-1-k}^{(r,s,m)}(t)+m(\ell-1)m_{\ell}^{(r,s,m)}(t).
	\end{align*}
	Proceeding in the same way, the second term may be written as
	\begin{align*}
		\sum_{i=1}^m\sum_{j\neq i}\frac{(x^i(t))^\ell+(x^i(t))^{\ell-1}x^j(t) }{x^i_t - x^j(t)} & =m^2\sum_{k=0}^{\ell-2}m_{k}^{(r,s,m)}(t)m_{\ell-1-k}^{(r,s,m)}(t)-m(\ell-1)m_{\ell-1}^{(r,s,m)}.
	\end{align*}
	Gathering all the terms, we get \eqref{ODEmoments}. The ODE satisfied by $m_1^{(r,s,m)}$ is derived exactly along the same lines noticing that the double sum is empty. 
\end{proof}
Theorem \ref{thm:maintheorem} is now proved as follows. Recall from \cite{MR2384475} that the moments
\begin{equation*}
	m^{(\lambda,\theta)}_{\ell}(t) := \frac{\tau(X_t^{\ell})}{\tau(P)}, \quad \ell \geq 1, \quad m_0(t) = 1,
\end{equation*}
of the free Jacobi process 
\begin{equation*}
	X_t = PU_tQU_t^{\star}P,
\end{equation*}
viewed as an operator in the compressed algebra $(P\mathscr{A}P, \tau/\tau(P))$, satisfy the differential system:
\begin{equation}\label{Jacmoments}
	\frac{d}{dt} m^{(\lambda,\theta)}_{\ell}(t) = -\ell m^{(\lambda,\theta)}_{\ell}(t) + \ell \theta m^{(\lambda,\theta)}_{\ell-1}(t) + \ell \lambda \theta \sum_{j=0}^{\ell-2}m^{(\lambda,\theta)}_{\ell-j-1}(t) [m^{(\lambda,\theta)}_{j}(t) - m^{(\lambda,\theta)}_{j+1}(t)],
\end{equation}
where $\tau(P) = \lambda \theta \in (0,1], \tau(Q) = \theta \in (0,1]$. 

We assume that for any $\ell \geq 1$, $m_{\ell}^{(r(m),s(m),m)}(0)$ converges as $m$ goes to infinity in the regime described in \eqref{eqn:asymptoticregime}. Then, for any time $t > 0$, 
 $$ 
 \tilde{m}_{\ell}^{(r(m),s(m),m)}(t):=m_{\ell}^{(r(m),s(m),m)}(t)\left(\frac{t}{d(m)}\right)
 $$ 
 converges to $m_{\ell}(t)$ as well. In fact let $L \geq 1$, then the $L-$ dimensional vector 
 $$(\tilde{m}_{\ell}^{(r(m),s(m),m)}(t), L\geq \ell \geq 1)$$ is the solution to a finite-dimensional differential system with matrix $M^{(m)}_L$, the coefficients are read off :
\begin{multline*}
	\frac{d}{dt}\tilde{m}_{\ell}^{(r,s,m)} (t)=\ell\frac{p-\ell+1}{d}  \tilde{m}_{\ell-1}^{(r,s,m)}-\ell\frac{d-\ell+1}{d} \tilde{m}_{\ell}^{(r,s,m)}
	+ \frac{m}{d}\ell\sum_{k=0}^{\ell-2}(\tilde{m}_{k}^{(r,s,m)}-\tilde{m}_{k+1}^{(r,s,m)})\tilde{m}_{\ell-1-k}^{(r,s,m)}.
\end{multline*}
where $ 1 \leq l \leq L$. As $m$ goes to infinity, $M^{(m)}_l$ converges in $\mathcal{M}_{L}^{(m)}(\mathbb{R})$ to the matrix $M_L$ of the differential system \eqref{ODEmoments} truncated for $1 \leq l \leq L$. Since the matrix exponential map is continuous and since the initial data $\tilde{m}_{\ell}^{(r(m),s(m),m)}(0)$ converges as well, then 
\begin{align*}
(\tilde{m}_{\ell}^{(r(m),s(m),m)}(t), L\geq \ell \geq 1) &= e^{tM_L}(\tilde{m}_{\ell}^{(r(m),s(m),m)}(0), L\geq \ell \geq 1) \\ 
&\hspace{0.5cm}\rightarrow e^{tM_L}(\tilde{m}_{\ell}^{(\lambda,\theta)}(0), L\geq \ell \geq 1) = (\tilde{m}_{\ell}^{(\lambda,\theta)}(t), L\geq \ell \geq 1).
\end{align*}
But weak convergence of compactly-supported probability measures is equivalent to the convergence of their moments, whence Theorem \ref{thm:maintheorem} follows. 
\begin{remark}
Notice that a key ingredient in the previous proof is the identification of the limit of the moment sequence $(\tilde{m}_{\ell}^{(r(m),s(m),m)(t)})_{\ell}$ as a moment sequence of a (unique) probability distribution.  \end{remark}

\subsection{Proof of Theorem \ref{thm:finitefreedifference}}    
\label{sec:finitefreediff}
We closely follow the line of arguments presented in \cite{arizmendi2024s} to prove Propositions 7.2, 7.3 and 7.8 there, based on the following representation of the free $S$ transform.
Let $v \in (0,1)$ and denote by ${\rm Dil}_v(\mu)$ the dilation by $v$ of the probability measure $\mu$. in particular, when the moments of $\mu$ exist, it holds that:
\begin{equation*}
         \int x^{\ell} {\rm Dil}_v(\mu)(dx)= v^{\ell}\int x^{\ell} \mu(dx).
\end{equation*}
    
Recall also the definition of \emph{fractional additive convolution power} \cite{shlyakhtenko2020fractional} in terms of the free cumulants. If $\mu$ is a compactly supported probability distribution, then $\mu^{\boxplus (1/v)}$ is the unique compactly supported probability measure whose free cumulants are the $(1/v)$-multiples of those of $\mu$. With these definitions \cite{belinschi2008remarkable}:
\begin{align}
\label{eqn:formulasandg}
S_{\mu}(-v) = -G_{{\rm Dil}_{v}(\mu^{\boxplus ({1}/{v})})}(0), \quad v \in (0,1).
\end{align}
where $G$ denotes the Cauchy transform.
Finally, recall the definition of the operators $\partial_{k|m}$ acting in monic polynomials with degree $m$:
\begin{align}
\label{eqn:formulafreecumulantsderivatives}
\partial^{k|m}p_m = \frac{1}{(m)_{m-k}} p^{(m-k)}, \quad 0 \leq k \leq m,
\end{align}
where $p^{(m-k)}$ is the $(m-k)$-th derivative of $p$.
\newline 

	\textbf{Case 1: We make the additional assumption that there exists $\varepsilon > 0$ and $\eta > 0$ such that all the roots of the polynomials $p_m$ are contained in $[\varepsilon,\eta)$}. \newline

    Under this assumption, the convergence $\mu\llbracket p_m\rrbracket\to\mu$ is equivalent to the convergence of all the moments of $\mu\llbracket p_m\rrbracket\to\mu$ toward those of $\mu$ and further to the convergence of the \emph{finite free cumulants of $\mu\llbracket p_m\rrbracket\to\mu$} toward those of $\mu$. This is out of the scope of the present article to give an account about finite free cumulants, we refere the reader to \cite{arizmendi2018cumulants}.

	We begin by relating the finite free $S$-transform to the Cauchy transform of the counting measure of the roots of derivatives of $p_m$, (see \cite{arizmendi2024s}, Lemma 6.11 2). Let $v \in (0,1)$.
	\begin{align*}
		m(\mathcal{S}^{(m)}_{p_m}(v+\frac{1}{m}) - \mathcal{S}_{p_m}^{(m)}(v)) &= m(\mathcal{S}^{(m)}_{p_m}(-\frac{\ceil{mv}+1}{m}) - \mathcal{S}^{(m)}_{p_m}(-\frac{\ceil{mv}}{m})) \\
		& = m(G_{\mu\llbracket \partial^{\ceil{mv}|} p_m \rrbracket}(0)-G_{\mu\llbracket \partial^{\ceil{mv}+1|m} p_m \rrbracket}(0))                          \\ 
		&=-G_{\nu^{(m)}_{p_m}(v)}(0),
	\end{align*}
	where $\nu^{(m)}_{p_m}(v) : = \mu\llbracket \partial^{\ceil{mv}|} p_m \rrbracket - \mu\llbracket \partial^{\ceil{mv}+1|m} p_m \rrbracket$.
	Observe that this measure has total variation equal to $2d$, it is therefore not uniformly bounded. Nevertheless, we will study the convergence of its Cauchy transform by using the moment method. By definition, the $\ell^{th}$ moment $M_\ell(\nu_{p_m}^{(m)}(v))$ of $\nu_{p_m}^{(m)}(v)$ is
    
	\begin{align}
        \label{eqn:momentspol}
		 M_\ell(\nu_{p_m}^{(m)}(v)) & = m(M_{\ell}( \partial^{\ceil{mv}+1|m}p_m) - M_{\ell}( \partial^{\ceil{mv}|m}p_m )) 	
    \end{align}
    
	Recall the following asymptotic expansion between the moments $M_{\ell}(p),\, \ell \geq 1$ and the \emph{finite free cumulants} $\kappa^{(m)}_{\ell}(p),\ell \geq 1$ of a monic polynomial $p$ of degree $m$ (see \cite{arizmendi2018cumulants}):

	\begin{align}
		\label{eqn:momentcumulants}
		M_{\ell}(p) =\sum_{\pi \in {\rm NC }(\ell)} \kappa^{(m)}_{\pi}(p) +\frac{2\ell}{m} \sum_{\substack{r+s = \ell \\ \pi \in S_{{\rm NC }}(r,s)}} r^{-1}s^{-1}\kappa^{(m)}_{\pi}(p)+ O(\frac{1}{m^2})
	\end{align}
	and the following simple formula relating the finite free cumulant of a polynomial $p$ and its derivatives $\partial^{k|m}p$ (see \cite{arizmendi2024s}) :
	\begin{align}
		\label{eqn:cumulants}
		\kappa^{(k)}_{\ell}(\partial^{k|m}(p)) =
		(\frac{k}{m})^{\ell-1}\kappa^{(m)}_{\ell}(p),\quad \ell \leq k.
	\end{align}
	Inserting \eqref{eqn:momentcumulants} in \eqref{eqn:momentspol} and using further equation \eqref{eqn:cumulants}, we obtain:

	\begin{align*}
		M_{\ell}(\nu^{(m)}_{p_m}(v)) & = m\bigg(\sum_{\pi \in {\rm NC }(\ell)} \bigg[(\frac{\ceil{mv}+1}{m})^{\ell-|\pi|}\kappa^{(m)}_{\pi}(p_m) - (\frac{\ceil{mv}}{m})^{\ell-|\pi|}\kappa_{\pi}^{(m)}(p_m)\bigg] \\
		                             & \hspace{2cm}+ 2\frac{\ell}{\ceil{mv}+1}\sum_{\substack{r+s = \ell                                                                                                      \\ \pi \in S_{{\rm NC }}(r,s)}} r^{-1}s^{-1}(\frac{\ceil{mv}+1}{m})^{\ell-|\pi|}\kappa^{(m)}_{\pi}(p_m) \\
		                             & \hspace{5cm}- 2\frac{\ell}{\ceil{mv}}\sum_{\substack{r+s = \ell                                                                                                        \\ \pi \in S_{{\rm NC }}(r,s)}} r^{-1}s^{-1}(\frac{\ceil{mv}}{m})^{\ell-|\pi|}\kappa_{\pi}^{(m)}(p_m) \bigg)  \\
		                             & + O(1/m)
	\end{align*}
	For the first term in the right-hand side of the equation above, with $\kappa_{\pi}(\mu),\pi \in {\rm {\rm NC }}(\ell)$ the partitioned free cumulants of $\mu$:
	\begin{align*}
		m\sum_{\pi \in {\rm NC }(\ell)} \bigg[(\frac{\ceil{mv}+1}{m})^{\ell-|\pi|} & \kappa^{(m)}_{\pi}(p) - (\frac{\ceil{mv}}{m})^{\ell-|\pi|}\kappa_{\pi}^{(m)}(p)\bigg]                                                    \\
		                                                                    & = m\sum_{\pi \in {\rm NC }(\ell)} \bigg[(\frac{\ceil{mv}}{m})^{\ell-|\pi|}((1+\frac{1}{\ceil{mv}})^{\ell-|\pi|}-1)\kappa^{(m)}_{\pi}(p_m)\bigg] \\
		                                                                    & = \sum_{\substack{\pi \in {\rm NC }(\ell)                                                                                                       \\ \pi \neq \hat{0}_{\ell}}} \bigg[(\ell-|\pi|)v^{\ell-|\pi|-1}\kappa_{\pi}(\mu)\bigg] + O(1/m)
	\end{align*}
	The sum of the second term and the third term are easily seen to contribute to a factor $O(\frac{1}{m})$. Hence, we obtain for the $\ell^{th}$ moment of $\nu_{p_m}^{(m)}(v)$ the following asymptotic expansion:
	\begin{align*}
		M_{\ell}(\nu_{p_m}^{(m)}(v)) & = \sum_{\substack{\pi \in {\rm NC }(\ell)                                              \\ \pi \neq \hat{0}_{\ell}}} \bigg[(\ell-|\pi|)v^{\ell-|\pi|-1}\kappa_{\pi}(\mu)\bigg]  + O({1}/{m}) \\
		                             & = \partial_{v} \sum_{\pi \in {\rm NC }(\ell)} v^{\ell-|\pi|}\kappa_{\pi}(\mu) + O(1/m) \\
		                             & = \partial_v m_{\ell}({\rm Dil}_v\mu^{\boxplus \frac{1}{v}}) + O(1/m)           \\
		                             & = m_{\ell}(\partial_v {\rm Dil}_v\mu^{\boxplus \frac{1}{v}}) + O(1/m)
	\end{align*}
    where in the last line $\partial_v {\rm Dil}_v\mu^{\boxplus \frac{1}{v}}$ is the weak limit of 
    $$
    \frac{1}{\varepsilon } ( {\rm Dil}_{v+\varepsilon}\mu^{\boxplus \frac{1}{v+\varepsilon}} - {\rm Dil}_{v}\mu^{\boxplus \frac{1}{v}})
    $$
    as $\varepsilon$ goes to zero.
    Under our assumptions, weak convergence is equivalent to the convergence of the moments (as an easy application of the Ascoli's theorem) and thus $\partial_v {\rm Dil}_v(\mu^{\boxplus 1/v})$ is the unique compactly supported signed measure such that 
    $$
    M_{\ell}(\partial_v {\rm Dil}_v\mu^{\boxplus \frac{1}{v}}) = \partial_v M_{\ell}({\rm Dil}_v\mu^{\boxplus \frac{1}{v}}) 
    $$
	Let us argue now that the Cauchy transform of $\nu_{p_m}^{(m)}(v)$ converges to the Cauchy transform of $\partial_v {\rm Dil}_v\mu^{\boxplus \frac{1}{v}}$ (analytical in the complementary of an interval $[\varepsilon(v),\eta(v)]$ with $\varepsilon(v) > 0$ (resp. $\eta(v)$) as a consequence of Lemma 5.6 and the proof of Proposition 5.5 in \cite{arizmendi2024s}; in the sequel $\varepsilon$ (resp. $\eta$) will denote the minimum (resp. maximum) between $\varepsilon(v)$ and $\varepsilon >0$ (resp. between $\eta$ and $\eta(v)$)). We observe first that:
	$$G_{\nu_p(v)} = \{G_{\nu^{(m)}_{p_m}(v)},\, m \geq 1\}$$ is a normal family of holomorphic functions on the domain $\Omega = [\varepsilon,1]^{c} \subset \hat{\mathbb{C}} = \mathbb{C} \cup \{+\infty\}$, by applying Montel's criterion (Fundamental Normality Test \cite{schiff1993normal}). Let $G\colon \Omega\to \hat{\mathbb{C}}$ be an accumulation point of $G_{\nu_p(v)}$. First, $G$ is not identically equal to $\infty$ since $G_{\nu^{(m)}_{p_m}(v)}(\infty)=0$ for all $m \geq 1$.
	Since
	$$
		\frac{1}{\ell+1!}\frac{m^{\ell+1}}{dz^{\ell+1}}\biggr|_{z=0}G_{\nu^{(m)}_{p_m}(v)}(1/z) = m_{\ell}(\nu_{p_m^{(m)}}) \rightarrow  m_{\ell}(\partial_v {\rm Dil}_v\mu^{\boxplus \frac{1}{v}}),
	$$
	we infer from the local uniform convergence and the Cauchy formula
	$$
		\frac{1}{(\ell+1)!}\frac{m^{\ell}}{dz^{\ell+1}}\biggr|_{z=0}G(1/z) = m_{\ell}(\partial_v {\rm Dil}_v\mu^{\boxplus \frac{1}{v}}).
	$$
	Hence, from unicity of analytic continuation $\{G_{\nu^{(m)}_{p_m}(v)},\, m \geq 1\}$ has a unique adherence point. Normality implies that $G_{\nu^{(m)}_{p_m}(v)}$ converges toward $G_{\partial_v {\rm Dil}_v\mu^{\boxplus \frac{1}{v}}}$ locally uniformly on $\Omega$.
    Since $0 \in \Omega$, we get, locally uniformly on $v$,
    \begin{multline*}
\lim_{m\rightarrow+\infty}\nabla^{(m)}\mathcal{S}^{(m)}_{p_m}(v) = \lim_{m\rightarrow+\infty}m(\mathcal{S}^{(m)}_{p_m}(v + \frac{1}{m}) - \mathcal{S}^{(m)}(v)) = -G_{\partial_v {\rm Dil}_v\mu^{\boxplus \frac{1}{v}}}(0)\\=-\partial_v G_{{\rm Dil}_v\mu^{\boxplus \frac{1}{v}}}(0) = -\partial_v {S}_{\mu}(-v).   
    \end{multline*}
    
	The result is proven.
	\newline

	\textbf{Case 2: The roots of the sequence $(p_m)_{g\geq 1}$ are contained in $(0,\eta]$.} \newline
    
    Let $v \in (0,1-\mu(\{0\}))$.
	From Theorem 1.2 in \cite{arizmendi2024s} (and the computations done in Case $1$), $\mu \llbracket \partial^{\ceil{mv} | m}\rrbracket$
	tends to ${\rm Dil}_{v}(\mu^{\boxplus \frac{1}{v}})$. Also, by Lemma 5.4 of \cite{arizmendi2024s}, for large enough $j$, $\mu \llbracket \partial^{\ceil{mv}|m} p_m\rrbracket $
	has support contained in some $[\varepsilon, \eta]$. We apply Case $1$ to infer 
	\begin{align}
	&\nabla^{(m)} \mathcal{S}^{(m)}_{p_m}(-v) =  \frac{m}{\ceil{mv}}\nabla^{(\ceil{mv})} \mathcal{S}^{(\ceil{mv})}_{\partial^{\ceil{mv}|m}p_m} (-1) 
	= \frac{1}{v}\partial_v  S_{{\rm Dil}_v \mu^{\boxplus \frac{1}{v}}}(-1). 
	\end{align}
	Given that $S_{{\rm Dil}_v \mu^{\boxplus \frac{1}{v}}}(z) = S_{\mu}(vz) $, we infer 
	$$
	\frac{1}{v}\partial_v  S_{{\rm Dil}_v \mu^{\boxplus \frac{1}{v}}}(-1) = -\partial_v  S_{\mu}(-v),
	$$
	and the result is proven.
	\newline

	\textbf{Case 3 : Unbounded support} 
    \newline

    This case is dealt with in exactly the same way as in Section 7.4 in \cite{arizmendi2024s}.

\subsection{Finite free \texorpdfstring{$T$}{T}-transform of the Frozen Jacobi process}
In this paragraph, we put together the results proved in the last two sections to obtain a differential-difference equation satisfied by the finite free
$T$ transform of $\chi^{(r,s,m)}_t$ of the Hermitian Jacobi process. Once obtained, we pass to the limit therein and derive again the PDE satisfied by the $T$ transform of the free Jacobi process. 

Let $e^{(r,s,m)}_n(t), 0 \leq n \leq m,$ be the elementary symmetric functions in the roots of $\chi^{(r,s,m)}(t)$: 
$$
\chi_t^{(r,s,m)}(x) = \sum_{n=0}^m (-1)^n e_n^{(r,s,m)}(t)x^{m-n}, \quad e_0^{(r,s,m)}(t) = 1.
$$
The following proposition provides a differential system satisfied by $e^{(r,s,m)}_n(t), 1 \leq n \leq m$. In this respect, the reader readily notices that it is simpler than the one 
derived in \cite{MR4092353}, Lemma 3.1., and corresponding to the Jacobi particle system in $[-1,1]$. 

\begin{proposition}
	\label{prop:diffequ}
	For any $1 \leq n \leq m$ and any time $t > 0$:
	$$
		\frac{d}{dt}e^{(r,s,m)}_{n}(t)  = -n(p+q-(n-1))e_{n}^{(r,s,m)}(t) 
		                                +(m-(n-1))(p-(n-1))e^{(r,s,m)}_{n-1}(t).
	$$
    \end{proposition}
\begin{proof}
	According to Proposition \ref{prop:heatequationjacobi}, for any time $t\geq 0$:
	\begin{equation}\label{edo-charac}
		\partial_t\chi_t^{(r,s,m)}(x) =-m(d-m+1)\chi_t^{(r,s,m)}(x) -\mathcal{L}_x^{(r,s)}[\chi_t^{(r,s,m)}](x).
	\end{equation}
	One readily infers from \eqref{edo-charac}:
	\begin{align*}
		\mathcal{L}_x^{(r,s)}\chi_t^{(r,s,m)}(x)= & x(1-x)\partial_{xx}^2 + [(r+1) - (r+s+2)x]\partial_x\chi_{t}^{(r,s,m)}                                                   \\
		=                                         & \sum_{k=0}^{m}(-1)^ke^{(r,s,m)}_{k}(t)\left(x(1-x)\partial_{xx}^2 x^{m-k}+ [(r+1) - (r+s+2)x]\partial_x x^{m-k}\right)
		\\=&\sum_{k=0}^{m}(-1)^ke^{(r,s,m)}_{k}(t)(m-k)(m-k-1) (x^{m-k-1}-x^{m-k})
		\\&+\sum_{k=0}^{m}(-1)^ke^{(r,s,m)}_{k}(t) (m-k)[(r+1)x^{m-k-1} - (r+s+2)x^{m-k}]
		\\=&-\sum_{k=0}^{m}(-1)^k(m-k)\{(m-k-1)+ (r+s+2)\}e^{(r,s,m)}_{k}(t) x^{m-k}
		\\&+ \sum_{k=0}^{m-1}(-1)^k(m-k)[r+1 +m-k-1]e^{(r,s,m)}_{k}(t)x^{m-k-1}
		\\=&-\sum_{k=0}^{m}(-1)^k(m-k)(m-k+ r+s+1)e^{(r,s,m)}_{k}(t) x^{m-k}
		\\& -\sum_{k=1}^{m}(-1)^k(m-k+1)[r +m-k+1]e^{(r,s,m)}_{k-1}(t)x^{m-k}.
	\end{align*}
	By equating coefficients on both sides of \eqref{edo-charac}, we obtain for any $k \ge 1$:
	
	\begin{align*}
		\frac{d}{dt}e^{(r,s,m)}_{k}(t)= & \left[-m(d-m+1) +(m-k)(m-k+ r+s+1)\right]e^{(r,s,m)}_{k}(t)
		\\&\hspace{5cm}+(m-k+1)(r +m-k+1)e^{(r,s,m)}_{k-1}(t)
		\\=&\left[-m(d-m+1) +(m-k)(p+q-m-k+1)\right]e^{(r,s,m)}_{k}(t)
		\\&\hspace{5cm}+(m-k+1)(p-k+1)e^{(r,s,m)}_{k-1}(t)
		\\=&-k(p+q-k+1)e^{(r,s,m)}_{k}(t)
		+(m-k+1)(p-k+1)e^{(r,s,m)}_{k-1}(t).
	\end{align*}
\end{proof}
Using this proposition, we obtain the following corollary (we recall the notation $\ceil{\cdot}$ for the ceiling function):
\begin{corollary} Let $z \in (0,1)$ and $t\geq 0$. Then
	\label{cor:odefiniteTtransform}
	\begin{multline*}
		\partial_t T_t^{(r,s,m)}(z)  = [2(m-\ceil{mz})-(p+q)]T_t^{(r,s,m)}(z) +p-2(m-\ceil{mz})                                                                                   \\
		                             +(p-m+\ceil{mz}+1))\frac{({m-\ceil{mz}})}{m}\frac{\nabla^{(m)}T_t^{(r,s,m)}(z)}{T_t^{(r,s,m)}(z+\frac{1}{m})} .
	\end{multline*}
\end{corollary}
\begin{proof}
	First of all, note that $x=0$ is not a root of $\chi_t^{(r,s,m)}$ since 
    \begin{equation*}
       \mathbb{E}\prod_{i=1}^m(\lambdaprbetai) > 0.
    \end{equation*}
    As a matter of fact, if $n = m-\ceil{mz}+1, z\in (0,1),$ then definition \eqref{Defi-T-Trans} reduces to:
\begin{align*}
		T_t^{(r,s,m)}(z) = \frac{n}{m-n+1}\frac{e^{(r,s,m)}_n(t)}{e^{(r,s,m)}_{n-1}(t)}.
	\end{align*}
	 We differentiate this expression in the time variable $t$ and get:
	\begin{align*}
		\partial_tT^{(r,s,m)}(z)=\frac{n}{(m-n+1)e^{(r,s,m)}_{n-1}(t) } \frac{de^{(r,s,m)}_n}{dt}(t)- \frac{T_t^{(r,s,m)}(z)}{e^{(r,s,m)}_{n-1}(t)} \frac{de^{(r,s,m)}_{n-1}}{dt}(t),
	\end{align*}
	where we used the fact that shifting $z \mapsto z + (1/m)$ amounts to shift $\ceil{mz} \mapsto \ceil{mz} + 1$ and in turn $n \mapsto n-1$. 
    Appealing to Proposition \ref{prop:diffequ}, we infer the following equation for the derivative of $T_t^{(r,s,m)}(z)$:
	\begin{align*}
		\partial_tT_t^{(r,s,m)}(z) & = (n(n-1)-n(p+q))T_t^{(r,s,m)}(z) +{n}(p-(n-1))                                                                                                                       \\& \hspace{1cm}-T_t^{(r,s,m)}(z)(n-1)((n-2)- (p+q)) -\frac{T_t^{(r,s,m)}(z)}{T_t^{(r,s,m)}(z+\frac{1}{m})}(p-(n-2))({n-1}).                    
        \end{align*}
	Writing
	$$
		\frac{T_t^{(r,s,m)}(z)}{T_t^{(r,s,m)}(z+\frac{1}{m})} = - \frac{\nabla^{(m)}T_t^{(r,s,m)}(z)}{mT_t^{(r,s,m)}(z+\frac{1}{m})} +1
	$$
	and keeping in mind the relation $n - 1 = m-\ceil{mz}$, we get the desired result after straightforward computations.
\end{proof}
The following Proposition \ref{thm:PDEfromfinitefree} derives the PDE satisfied by the $T$ transform of the free Jacobi
process as the limit of the differential-difference equation stated in Proposition \ref{cor:odefiniteTtransform}.
\begin{proposition}
	\label{thm:PDEfromfinitefree}
	In the asymptotic regime \eqref{eqn:asymptoticregime},
	\begin{itemize}
		\item For each time $t\geq 0$, the $T$ transform $T^{(r,s,m)}_{t/d(m)}(\cdot)$ converges locally uniformly on $(-1,0)$.
		\item Let $\mathcal{T}^{(\lambda,\theta)}\colon \mathbb{R}\times (0,1) $ be the limit, then $\mathcal{T}^{(\lambda,\theta)}$ is continuously differentiable and it satisfies:
	\end{itemize}
	\begin{align*}
		&\partial_t \mathcal{T}_t^{(\lambda,\theta)}(z) = (2(1-z)\lambda\theta -1)\mathcal{T}_t^{(\lambda,\theta)}(z) + \theta(1-2\lambda (1-z)) \\
		&\hspace{5cm}+ \theta(1-\lambda{(1-z)})(1-z){\partial_{z}\log \mathcal{T}_t^{(\lambda,\theta)}(z)}, \\
		&\mathcal{T}_0^{(\lambda,\theta)}(z) = 1
	\end{align*}
\end{proposition}
\begin{proof}
We already know that for each time $t \geq 0$, the finite free $T$ transform $T^{(r,s,m)}_{t/d(m)}$ converges. By virtue of Corollary \ref{cor:odefiniteTtransform}, the only point missing is the local uniform convergence (in time $t$, as $m$ goes to infinity) of $T^{(r,s,m)}_{t/d(m)},~t\geq 0$.  

We show that the root distribution $\mu^{(r,s,m)}_{t/d(m)}$ converges weakly, locally uniformly in time, which means that for any smooth and bunded function $f: [0,1]\to\mathbb{R}$, 
$
\mu_{t/d(m)}^{(r,s,m)}(f),\, t\geq 0,
$
converges locally uniformly in time. The point-wise convergence in time of the $T$-transform implies that for every $t \geq 0$, $\mu_{t/d(m)}^{(r,s,m)}(f)$ converges (in the regime \ref{eqn:asymptoticregime}). Let us compute the derivative in time:
\begin{multline*}
		\partial_t\mu^{(r,s,m)}_{t/d(m)}(f) = \frac{p}{d(m)} + \frac{p+q}{d(m)}\left\{\frac{1}{m}\sum_{j}x_jf'(x_j) \right\} \\
		+ \frac{1}{2md(m)} \sum_{j\neq k} \frac{f(x_j)-f(x_k)}{x_j-x_k}(x_j(1-x_k)+x_k(1-x_j)).
\end{multline*}
Furthermore, we can bound this derivative uniformly on $m$ by
\begin{align*}
		|\partial_t\mu^{(r,s,m)}_{t/d(m)}(f)| & \leq \frac{p}{d(m)} + \frac{p+q}{d(m)} \|f'\|_{\infty} + \frac{1}{md(m)} \sum_{k,j} \frac{|f(x_j)-f(x_k)|}{|x_j-x_k|} \\
		                               & \leq \frac{p}{d(m)} + \frac{p+q}{d(m)} \|f'\|_{\infty} + \frac{m(m-1)}{md(m)}\|f''\|_{\infty}.
\end{align*}
We infer that, for any $f$ smooth,
$
t\mapsto \mu_{t/d(m)}^{(r,s,m)}(f)
$
is equicontinuous and the modulus of continuity is uniform in $m$. Therefore, pointwise convergence in time and continuity of the limit imply the local uniform convergence of $\mu^{(r,s,m)}_{t/d(m)}(f)$ as $m$ goes to infinity. This holds, in particular, for the moments of $\mu^{(r,s,m)}_{t/d(m)}$ and thus for all polynomial functions of the moments $\mu^{(r,s,m)}_{t/d(m)}$ ; the $e^{(r,s,m)}_n(t)$; 
$m \geq n \geq 0$. 

We infer that all finite free cumulants of $\chi^{(r,s,m)}_t$ converge locally uniformly in time. 

From the formula recalled in equation \eqref{eqn:formulafreecumulantsderivatives}, we also deduce that all the finite free cumulants of the iterated derivatives of $\chi_{t/d(m)}^{(r,s,m)}$ converge locally uniformly in time (see equations \eqref{eqn:cumulants}) and therefore the weak convergence, locally uniformly in time, of their root distributions.

The formula \eqref{eqn:formulasandg} will imply local uniform convergence in time of the finite free $S$ and $T$ - transform of $\chi_{t/d(m)}^{(r,s,m)}$, provided that we can bound away from zero the roots of the iterated derivatives of $\chi_{t/d(m)}^{(r,s,m)}$, locally uniformly in time (since $S$ and $T$ are rational expressions in the $e_n^{(r,s,m)}(t)$). 
This is done using Lemma 5.5 in \cite{arizmendi2024s}, since the cutoff $a$ the author chooses can be chosen uniformly in time, locally for the distributions $\nu_t^{(\lambda,\theta)}$ : if $t_0 \in \mathbb{R}$, we can choose $a$ such that the c.d.f $F_{t_0}$ of $\nu^{(\lambda,\theta)}_{t_0}$ is continuous at $a$ and $F_{t_0}(a) < 1-u$, under the notations used by the authors in \cite{arizmendi2024s}. We can pick a smooth function $f:(0,1)\to\mathbb{R}$ that is equal to $1$ on $(a,1)$ and equal to $0$ on $(0,a-\eta)$ such that 
\begin{equation}
\label{eqn:elementaryinequality}
F_{t_0}(a) < \nu^{(\lambda,\theta)}_{t_0}(f) < 1-u.
\end{equation}
where the first inequality follows from $f > 1_{[a,1)}$. Since $\nu_{t}^{\lambda,\theta}(f)$ is continuous at $t_0$, the inequality \ref{eqn:elementaryinequality}
holds for $t$ in the vicinity of $t_0$.

Together we the second Dini's Theorem, since $T_{t/d(m)}^{(r,s,m)}(\cdot)$ is a non-decreasing step-function, we deduce that the right-hand side of the differential-difference equation of Corollary \ref{cor:odefiniteTtransform}, if normalized by the factor $d(m)$ and at times $t/d(m)$, converges locally uniformly in time.

The limit $T^{(\lambda,\theta)}(z)$ is therefore continuously differentiable in time and for each time $t\geq 0$ and $z \in (-1,0)$:
$$
\partial_t T^{(r,s,m)}_{t/d(m)}(z) \rightarrow T^{(\lambda,\theta)}_t(z)
$$
and we can pass to the limit on both sides of the integro-differential equation of Corollary \ref{cor:odefiniteTtransform} to obtain the second assertion of the theorem. 
\end{proof}
\printbibliography
\end{document}